\documentclass[10pt,reqno]{amsart}
   \usepackage[centertags]{amsmath}
   \usepackage{amsfonts}
   \usepackage{amsmath}
   \usepackage{amssymb}
   \usepackage{amsthm}
   \usepackage{newlfont}
   \usepackage[latin1]{inputenc}
\usepackage{multirow, bigdelim}

\usepackage{epsfig}
\usepackage{pst-all}
\usepackage{pstricks}
\usepackage{pgf}
\usepackage{mathrsfs}
\usepackage{hyperref}
\usepackage{bm}

\usepackage{graphicx}

\usepackage{bm}

\usepackage{mathtools}
\usepackage{ifmtarg}

\makeatletter

    \newcommand\contFrac{\@ifstar{\@contFracStar}{\@contFracNoStar}}

    \def\singleContFrac#1#2{%
        \begin{array}{@{}c@{}}%
            \multicolumn{1}{c|}{#1}%
            \\%
            \hline%
            \multicolumn{1}{|c}{#2}%
        \end{array}%
    }

    \def\@contFracNoStar#1{%
        \mathchoice{
            \@contFracNoStarDisplay@#1//\@nil%
        }{
            \@contFracNoStarInline@#1//\@nil%
        }{
            \@contFracNoStarInline@#1//\@nil%
        }{
            \@contFracNoStarInline@#1//\@nil%
        }%
    }

    \def\@contFracNoStarDisplay@#1//#2\@nil{%
        \@ifmtarg{#2}{%
            #1%
        }{%
            #1+\cfrac{1}{\@contFracNoStarDisplay@#2\@nil}%
        }%
    }

        \def\@contFracNoStarInline@#1//#2\@nil{%
            \@ifmtarg{#2}{%
                #1%
            }{%
                #1 \@@contFracNoStarInline@@#2\@nil%
            }%
        }
        \def\@@contFracNoStarInline@@#1//#2\@nil{%
            \@ifmtarg{#2}{%
                + \singleContFrac{1}{#1}%
            }{%
                + \singleContFrac{1}{#1} \@@contFracNoStarInline@@#2\@nil%
            }%
        }

    \def\@contFracStar#1{%
        \mathchoice{
            \@contFracStarDisplay@#1////\@nil%
        }{
            \@contFracStarInline@#1//\@nil%
        }{
            \@contFracStarInline@#1//\@nil%
        }{
            \@contFracStarInline@#1//\@nil%
        }%
    }

    \def\@contFracStarDisplay@#1//#2//#3\@nil{%
        \@ifmtarg{#2}{%
            #1%
        }{%
            #1 + \cfrac{#2}{\@contFracStarDisplay@#3\@nil}%
        }%
    }

        \def\@contFracStarInline@#1//#2\@nil{%
            \@ifmtarg{#2}{%
                #1%
            }{%
                #1 \@@contFracStarInline@@#2\@nil%
            }%
        }
        \def\@@contFracStarInline@@#1//#2//#3\@nil{%
            \@ifmtarg{#3}{%
                - \singleContFrac{#1}{#2}%
            }{%
                - \singleContFrac{#1}{#2} \@@contFracStarInline@@#3\@nil%
            }%
        }
\makeatother

\allowdisplaybreaks[3]

\theoremstyle{plain}
\newtheorem{theorem}{Theorem}[section]

\theoremstyle{definition}

\theoremstyle{remark}
\newtheorem{remark}[theorem]{Remark}
\newtheorem*{remark*}{Remark}

\numberwithin{equation}{section}

\title[Stochastic Darboux for QBD processes and urn models]{Stochastic Darboux transformations for \\quasi-birth-and-death processes and urn models}

\author{F. Alberto Gr\"unbaum}
\address{F. Alberto Gr\"unbaum\\
Department of Mathematics. University of California, Berkeley. Berkeley, CA 94720  U.S.A.}
\email{grunbaum@math.berkeley.edu}
\date{\today}
\thanks{
}

\author{Manuel D. de la Iglesia}
\address{Manuel D. de la Iglesia\\
Instituto de Matem\'aticas, Universidad Nacional Aut\'onoma de M\'exico, Circuito Exterior, C.U., 04510, Ciudad de México, México.}
\email{mdi29@im.unam.mx}
\thanks{The work of both authors is supported by UC MEXUS-CONACYT grant CN-16-84. The work of the first author is also supported by PAPIIT-DGAPA-UNAM grant IA102617 (México), MTM2015-65888-C4-1-P (Ministerio de Economía y Competitividad, Spain) and FQM-262, FQM-7276 (Junta de Andaluc\'ia).}

\date{\today}



\begin{document}

\maketitle

\begin{abstract}
We consider stochastic UL and LU block factorizations of the one-step transition probability matrix for a discrete-time quasi-birth-and-death process, namely a stochastic  block tridiagonal matrix. The simpler case of random walks with only nearest neighbors transitions gives a unique LU factorization and a one-parameter family of factorizations in the UL case. The block structure considered here yields many more possible factorizations resulting in a much enlarged class of potential applications. By reversing the order of the factors (also known as a Darboux transformation) we get new families of quasi-birth-and-death processes where it is possible to identify the matrix-valued spectral measures in terms of a  Geronimus (UL) or a Christoffel (LU) transformation of the original one. We apply our results to one example going with matrix-valued Jacobi polynomials arising in group representation theory. We also give urn models for some particular cases.
\end{abstract}

\section{Introduction}

Among the class of Markov chains there is one set  that can be analyzed by so-called ``spectral methods'', namely random walks (discrete-time) and birth-and-death processes (continuous time). They go with a one-step tridiagonal matrix and this naturally leads to a self-adjoint operator in certain Hilbert space. Starting with \cite{KMc2,KMc3,KMc6}, and using earlier ideas of W. Feller and H. P. McKean in the case of diffusion processes, there is a vast literature on this subject, which relies on the rich theory of orthogonal polynomials. A historical overview of this material can be seen, for instance, in \cite{ILMV}. Eventually this approach was extended to cover so-called quasi-birth-and-death (QBD) processes by exploiting matrix-valued orthogonal polynomials, a notion due to M. G. Kre$\breve{\mbox{{\i}}}$n, see \cite{K2,K1} for this, and \cite{DRSZ,G2} for its use in the study of QBD processes. Here the tridiagonal matrix is replaced by a block tridiagonal one. For a general reference about QBD processes see \cite{LaR, Neu}. The spectral methods work well when one has knowledge of the orthogonal polynomials and the spectral measure associated with the one-step transition probability matrix. Needless to say, this is a limitation on the wide practical use of the method, although many interesting general results are available.

\smallskip

On the other hand, given either a tridiagonal or a block tridiagonal matrix, assumed here to be stochastic, a natural problem is to explicitly construct some simple probabilistic model (such as an urn model) whose one-step transition probability matrix coincides with the given one. This is not as simple as it may sound: to give an instance of this consider the  case of the matrices corresponding to the Krawtchouk and the Hahn orthogonal polynomials. In these two cases the urn model predates the consideration of the
spectral problem for the tridiagonal matrices by a long stretch, see for instance \cite{G2} or \cite{F}, p. 378, where one sees the connection with work of P. and T. Ehrenfest in 1907 (see \cite{E}), as well as D. Bernoulli (1769) and S. Laplace (1812). Both of these cases are very special cases of the so-called Askey-Wilson tridiagonal matrices, see \cite{AW}, for which nobody has been able to find a nice urn model. For another well known example of a tridiagonal matrix, namely the one associated to the Jacobi polynomials, a first and rather contrived urn model was only given, as far as we know, in \cite{G4} (see also \cite{GdI3} for a different urn model). For a matrix-valued version of the Jacobi polynomials a pair of explicit models was given in \cite{GPT3}.

\smallskip

Keeping in mind the two points raised above, we are in a position to describe the purpose of this paper. Extending our previous work for random walks in \cite{GdI3}, we consider the block tridiagonal transition probability matrix of a discrete-time QBD process and perform UL and LU stochastic block factorizations. Unlike the case of random walks, where the UL factorization comes with exactly one free parameter, now, for the stochastic block factorization, there may be many degrees of freedom, as we will see in Section \ref{SEC2}.  The same applies for the LU factorization. The main motivation of this factorization is to analyze the urn model associated with the QBD process in terms of two unrelated and simpler urn models and to combine them to obtain a simpler description of the original QBD process.

\smallskip

Once we are able to perform UL and LU stochastic block factorizations of the block tridiagonal transition probability matrix of a discrete-time QBD process, we will give a general way to produce a family of new ones, performing what is called a discrete Darboux transformation by reversing the order of multiplication (see Section \ref{SEC3}). We will also give a way to relate the original and the new spectral ingredients, i.e. the matrix-valued orthogonal polynomials and the matrix-valued spectral measure.

\smallskip

We apply our results in Section \ref{SEC4}  to study one example of Jacobi type coming from group representation theory, introduced for the first time in \cite{GPT1} (see also \cite{GPT3}). We focus on the $2\times2$ case and study two particular situations, where we can illustrate the features that arise in the case of a general QBD process. Finally, in Section \ref{SEC5}, we start from a special case of the urn model described in \cite{GPT3} and find a different urn model as an application of the method of the stochastic block factorization.

\section{Stochastic LU and UL block factorization}\label{SEC2}

Let $P$ be the one-step transition probability matrix of a discrete-time quasi-birth-and-death (QBD) process with state space $\mathbb{Z}_{\geq0}\times\{1,2,\ldots,d\}$, $d\geq1$, given by
\begin{equation}
P=\begin{pmatrix}\label{PP}
B_0&A_0&0\\
C_1&B_1&A_1&\\
0&C_2&B_2&A_2&\\
&&\ddots&\ddots&\ddots
\end{pmatrix}.
\end{equation}
Here $A_n$, $B_n$ and $C_n$ are $d\times d$ matrices. We will assume for simplicity that the matrices $A_n$ and $C_n$ are nonsingular. The symbol $0$ as well as all the unfilled entries in \eqref{PP} stand for the $d\times d$ block zero matrix. When $d=1$ we recover the classical random walk with state space $\mathbb{Z}_{\geq0}$.

Let us denote by $e_j$ the $j$-th canonical $d$-dimensional vector and $\bm e_d$ the $d$-dimensional vector with all components equal to 1, i.e.
$$
\bm e_d=(1,1,\ldots,1)^T.
$$
Since $P$ is a stochastic matrix, we have nonnegative (scalar) entries, i.e.
$$
e_i^TA_ne_j\geq0,\quad e_i^TB_ne_j\geq0,\quad e_i^TC_{n+1}e_j\geq0,\quad i,j=1,\ldots,d, \quad n\geq0,
$$
and all (scalar) rows add up to one, i.e.
$$ 
(B_0+A_0)\bm e_d=\bm e_d,\quad (C_n+B_n+A_n)\bm e_d=\bm e_d,\quad n\geq1.
$$
Observe that all block entries of $P$ are \emph{semi-stochastic} $d\times d$ matrices, i.e. all entries are nonnegative and $A_n\bm e_d\leq\bm e_d$, $B_n\bm e_d\leq\bm e_d$ and $C_{n+1}\bm e_d\leq\bm e_d$ for $n\geq0$ (component wise).

\vspace{1.6cm}

A diagram of the transitions between states looks as follows (for $d=2$)

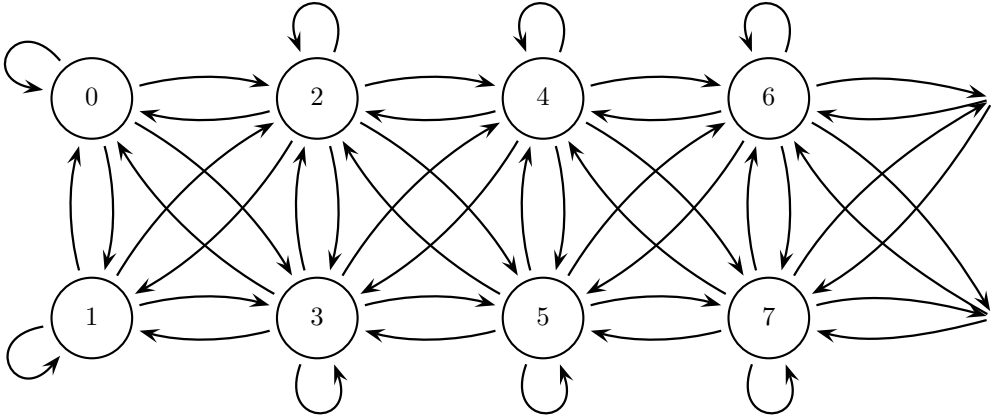
\vspace{1.0cm}

$$\begin{psmatrix}[rowsep=2.5cm,colsep=3cm]
  \cnode{.55}{0}& \cnode{.55}{2} & \cnode{.55}{4} & \cnode{.55}{6} & \pnode{8} \\
  \cnode{.55}{1} & \cnode{.55}{3} & \cnode{.55}{5} & \cnode{.55}{7} & \pnode{9} \\
\psset{nodesep=3pt,arcangle=15,labelsep=2ex,linewidth=0.3mm,arrows=->,arrowsize=1mm
3} \nccurve[angleA=130,angleB=170,ncurv=4]{0}{0}
\nccurve[angleA=190,angleB=230,ncurv=4]{1}{1}
\nccurve[angleA=70,angleB=110,ncurv=4]{2}{2}
\nccurve[angleA=70,angleB=110,ncurv=4]{4}{4}
\nccurve[angleA=70,angleB=110,ncurv=4]{6}{6}
\nccurve[angleA=250,angleB=290,ncurv=4]{3}{3}
\nccurve[angleA=250,angleB=290,ncurv=4]{5}{5}
\nccurve[angleA=250,angleB=290,ncurv=4]{7}{7} \ncarc{0}{2}
\ncarc{2}{0} \ncarc{2}{4} \ncarc{4}{2} \ncarc{4}{6} \ncarc{6}{4}
\ncarc{6}{8} \ncarc{8}{6} \ncarc{0}{1} \ncarc{1}{0} \ncarc{2}{1}
\ncarc{1}{2} \ncarc{1}{3} \ncarc{3}{1} \ncarc{2}{3} \ncarc{3}{2}
\ncarc{4}{3} \ncarc{3}{4}\ncarc{3}{5} \ncarc{5}{3}\ncarc{5}{4}
\ncarc{4}{5}\ncarc{5}{6} \ncarc{6}{5}\ncarc{5}{7}
\ncarc{7}{5}\ncarc{6}{7} \ncarc{7}{6}\ncarc{7}{8}
\ncarc{8}{7}\ncarc{7}{9} \ncarc{9}{7} \ncarc{0}{3}\ncarc{3}{0}
\ncarc{2}{5}\ncarc{5}{2}\ncarc{4}{7}\ncarc{7}{4}\ncarc{6}{9}\ncarc{9}{6}
\psset{labelsep=-4.25ex}\nput{90}{0}{0}
\psset{labelsep=-4.25ex}\nput{90}{2}{2}
\psset{labelsep=-4.25ex}\nput{90}{4}{4}
\psset{labelsep=-4.25ex}\nput{90}{6}{6}
\psset{labelsep=-4.25ex}\nput{90}{1}{1}
\psset{labelsep=-4.25ex}\nput{90}{3}{3}
\psset{labelsep=-4.25ex}\nput{90}{5}{5}
\psset{labelsep=-4.25ex}\nput{90}{7}{7}
\end{psmatrix}
$$

\vspace{-1.5cm}

It will turn out to be useful to perform a UL block factorization of the matrix $P$ in the following way
\begin{equation}\label{Pxy}
P=\begin{pmatrix}
B_0&A_0&\\
C_1&B_1&A_1&\\
&\ddots&\ddots&\ddots
\end{pmatrix}=\begin{pmatrix}
Y_0&X_0&\\
0&Y_1&X_1&\\
&\ddots&\ddots&\ddots
\end{pmatrix}\begin{pmatrix}
S_0&0&\\
R_1&S_1&0&\\
&\ddots&\ddots&\ddots
\end{pmatrix}=P_UP_L,
\end{equation}
with the condition that $P_U$ and $P_L$ \emph{are also stochastic matrices}, i.e. all (scalar) entries are nonnegative and 
\begin{equation}\label{stXYRS} 
(X_n+Y_n)\bm e_d=\bm e_d,\quad n\geq0,\quad S_0\bm e_d=\bm e_d,\quad (R_n+S_n)\bm e_d=\bm e_d,\quad n\geq1.
\end{equation}
Observe that $S_0$ must be stochastic, while the rest of the block entries of $P_U$ and $P_L$ must be semi-stochastic. A direct computation shows that
\begin{align}
\nonumber A_n&=X_nS_{n+1},\quad n\geq0,\\
\label{ULrel} B_n&=X_nR_{n+1}+Y_nS_n,\quad n\geq0,\\
\nonumber C_n&=Y_nR_n\quad n\geq1.
\end{align}
Since $A_n$ and $C_n$ are nonsingular then $X_n, Y_n,R_{n+1}$ and $S_{n+1}, n\geq0$, are nonsingular matrices. $S_0$ may or may not be nonsingular. This factorization, just as in the scalar situation (see \cite{GdI3}), simplifies the interpretation of the original QBD process $P$, expressing it as the composition of two simpler QBD processes, $P_U$ and $P_L$.

One could have performed the factorization the other way around in the form
\begin{equation}\label{PxyLU}
P=\begin{pmatrix}
B_0&A_0&\\
C_1&B_1&A_1&\\
&\ddots&\ddots&\ddots
\end{pmatrix}=\begin{pmatrix}
\widetilde S_0&0&\\
\widetilde R_1&\widetilde S_1&0&\\
&\ddots&\ddots&\ddots
\end{pmatrix}\begin{pmatrix}
\widetilde Y_0&\widetilde X_0&\\
0&\widetilde Y_1&\widetilde X_1&\\
&\ddots&\ddots&\ddots
\end{pmatrix}=\widetilde P_L\widetilde P_U,
\end{equation}
in which case we have a LU block factorization with relations
\begin{align*}
A_n&=\widetilde S_{n}\widetilde X_n,\quad n\geq0,\\
B_n&=\widetilde R_{n}\widetilde X_{n-1}+\widetilde S_n\widetilde Y_n,\quad n\geq0,\\
C_n&=\widetilde R_n\widetilde Y_{n-1}\quad n\geq1.
\end{align*}
As we will see the LU block factorization will have \emph{fewer} degrees of freedom than the UL block factorization.

\medskip

Let us focus first on the UL block factorization \eqref{Pxy}. We remark here an important difference with respect to the scalar situation (see \cite{GdI3}). In the scalar situation the UL factorization has exactly one free parameter $y_0$ (defined below), while in the LU factorization case the factorization is unique. This is \emph{not the case} for the UL and LU block factorizations, where there may be many degrees of freedom. For instance, in the scalar situation, where we use lower case symbols, one could compute $x_n$ and $r_{n+1}$ in terms of $y_n$ and $s_{n+1},$ respectively, for $n\geq0$, since both factors $P_U$ and $P_L$ must be stochastic matrices. Then $x_n=1-y_n, r_{n+1}=1-s_{n+1}, n\geq0,$ and $s_0=1$, i.e. all entries of $P_U$ and $P_L$ can be computed in terms of only one free parameter, $y_0$ (see Lemma 2.2 of \cite{GdI3}).

The stochasticity conditions on $P_U$ and $P_L$ gives the relations \eqref{stXYRS}, but we notice that it is not possible to compute, for instance, all entries of $S_0$ by having only the information that $S_0\bm e_d=\bm e_d$. The same is true for the rest of coefficients, i.e. we can not compute all entries of $X_n$ in terms of $Y_n$ from the information that $(X_n+Y_n)\bm e_d=\bm e_d$ (same for $R_n$ and $S_n$). 

One way of computing the block entries $X_n, Y_n, R_{n}, S_n$ comes directly from \eqref{ULrel}. From the first and third relation we can compute $S_n, R_n, n\geq1$ in terms of $X_n, n\geq0$ and $Y_n, n\geq1$, respectively. The second relation gives then
$$
Y_{n+1}=C_{n+1}\left(B_n-Y_nX_{n-1}^{-1}A_{n-1}\right)^{-1}X_n,\quad n\geq1,\quad\mbox{and}\quad Y_{1}=C_{1}\left(B_0-Y_0S_0\right)^{-1}X_0.
$$

Therefore, all coefficients can be computed in terms of $Y_0, S_0$ and $X_n,n\geq0$. 
We will see below that all these inverses are well defined as long as 
certain invertibility conditions are satisfied. Apart from this we have to impose certain positivity conditions on all entries of $X_n, Y_n, R_{n}, S_n$. This gives an infinite number of free parameters in general and it is very difficult to pick a ``natural one'' among the possible solutions.


\medskip

We propose now one way of computing the block entries $X_n, Y_n, S_n, R_n$ using what corresponds to the ``monic'' version of $P$ in connection with matrix-valued polynomials, see \cite{K2,K1}. This method was used, for a different purpose, in \cite{G3}. Let us call
\begin{equation}\label{LLn}
L_n=(A_0\cdots A_{n-1})^{-1},\quad n\geq1,\quad L_0=I.
\end{equation}
Observe that all $L_n$ are invertible, by assumption, and $L_nL_{n+1}^{-1}=A_n$. Then we have
\begin{equation}\label{PLJ}
P=LJL^{-1},
\end{equation}
where
$$
L=\begin{pmatrix}
L_0&\\
&L_1&\\
&&\ddots
\end{pmatrix},\quad J=\begin{pmatrix}
\widehat{B}_0&I&\\
\widehat{C}_1&\widehat{B}_1&I&\\
&\ddots&\ddots&\ddots
\end{pmatrix}.
$$
The ``monic'' block entries $\widehat{B}_n$ and $\widehat{C}_n$ are related to the old ones by the relations
\begin{equation*}\label{Cow}
\widehat{B}_n=L_n^{-1}B_nL_n,\quad n\geq0,\qquad \widehat{C}_n=L_n^{-1}C_nL_{n-1}=L_n^{-1}C_nA_{n-1}L_n,\quad n\geq1.
\end{equation*}
Consider now the UL block factorization of the ``monic'' operator $J$ in the following way
\begin{equation}\label{Jab}
J=\begin{pmatrix}
\widehat{B}_0&I&\\
\widehat{C}_1&\widehat{B}_1&I&\\
&\ddots&\ddots&\ddots
\end{pmatrix}=\begin{pmatrix}
\alpha_0&I&\\
0&\alpha_1&I&\\
&\ddots&\ddots&\ddots
\end{pmatrix}\begin{pmatrix}
I&0&\\
\beta_1&I&0&\\
&\ddots&\ddots&\ddots
\end{pmatrix}=\bm\alpha\bm\beta.
\end{equation}
Then we have
\begin{align}
\label{alfbet}\widehat{B}_{n}&=\beta_{n+1}+\alpha_{n},\quad n\geq0,\\
\nonumber\widehat{C}_{n}&=\alpha_n\beta_n,\quad n\geq1.
\end{align}
These relations give a direct way to compute $\alpha_n$ and $\beta_n$ in terms of $\widehat{B}_{n}$, $\widehat{C}_{n}$ and a \emph{free matrix-valued parameter} $\alpha_0$. Indeed, the first terms are given recursively by
\begin{align*}
\beta_1&=\widehat{B}_{0}-\alpha_0,& \alpha_1&=\widehat{C}_{1}(\widehat{B}_{0}-\alpha_0)^{-1}\\
\beta_2&=\widehat{B}_{1}-\widehat{C}_{1}(\widehat{B}_{0}-\alpha_0)^{-1},&\alpha_2&=\widehat{C}_{2}\left[\widehat{B}_{1}-\widehat{C}_{1}(\widehat{B}_{0}-\alpha_0)^{-1}\right]^{-1},\quad\mbox{etc.}
\end{align*}
Observe from \eqref{alfbet}, and since $\widehat{C}_n$ is a nonsingular matrix, that $\alpha_n,\beta_n$ must be nonsingular matrices as well. This means that the free parameter $\alpha_0$ is subject to certain invertibility conditions, i.e. $\widehat{B}_{0}-\alpha_0, \widehat{B}_{1}-\widehat{C}_{1}(\widehat{B}_{0}-\alpha_0)^{-1}$, etc, must be nonsingular matrices.

Substituting \eqref{Jab} into \eqref{PLJ} leads to
$$
P=\left[L\bm\alpha\right]\left[\bm\beta L^{-1}\right],
$$
which is a UL factorization of $P$. In general, this factorization will not give stochastic factors as in \eqref{Pxy}. To guarantee this we will need to introduce below more degrees of freedom while keeping the UL structure.  Let us denote by  $T$ the block diagonal invertible matrix
$$
T=\begin{pmatrix}
\tau_0&\\
&\tau_1&\\
&&\ddots
\end{pmatrix}.
$$
Then we can write $P$ as
\begin{equation*}
P=\left[L\bm\alpha T\right]\left[T^{-1}\bm\beta L^{-1}\right]=P_UP_L.
\end{equation*}
Identifying block entries with \eqref{Pxy} we get
\begin{align}
\label{xyrsg}X_n&=L_n\tau_{n+1},\quad Y_n=L_n\alpha_n\tau_n,\quad n\geq0,\\
\nonumber S_n&=\tau_n^{-1}L_n^{-1},\quad R_{n+1}=\tau_{n+1}^{-1}\beta_{n+1}L_{n}^{-1},\quad n\geq0.
\end{align}
Relations \eqref{stXYRS} give
\begin{align}
\label{eq1} L_n(\alpha_n\tau_n+\tau_{n+1})\bm e_d&=\bm e_d,\quad n\geq0\\
\label{eq2} (\beta_{n+1}L_{n}^{-1}+L_{n+1}^{-1})\bm e_d&=\tau_{n+1}\bm e_d,\quad n\geq0,\\
\label{eq3} \tau_0^{-1}\bm e_d&=\bm e_d.
\end{align}
If \eqref{eq2} and \eqref{eq3} hold, then \eqref{eq1} must hold as well. Indeed, for $n=0$, we have
$$
(\alpha_0\tau_0+\tau_{1})\bm e_d=(\alpha_0+\beta_1+L_1^{-1})\bm e_d=(B_0+A_0)\bm e_d=\bm e_d,
$$
while for $n\geq1$, we have
\begin{align*}
L_n(\alpha_n\tau_n+\tau_{n+1})\bm e_d&=\left[L_n\alpha_n(\beta_{n}L_{n-1}^{-1}+L_{n}^{-1})+L_n(\beta_{n+1}L_{n}^{-1}+L_{n+1}^{-1})\right]\bm e_d\\
&=\left[L_n\alpha_n\beta_{n}L_{n-1}^{-1}+L_n\alpha_nL_{n}^{-1}+L_n\beta_{n+1}L_{n}^{-1}+L_nL_{n+1}^{-1}\right]\bm e_d\\
&=\left[C_n+B_n+A_n\right]\bm e_d=\bm e_d.
\end{align*}
Observe that \eqref{eq2} gives one way to compute $\tau_n$ which may allow for several degrees of freedom as well. Computing $\tau_n$ using \eqref{eq2} and \eqref{eq3} does not yet guarantee that $P_U$ and $P_L$ are going to be stochastic. We must still verify the positivity of the entries.

Using \eqref{xyrsg} and the explicit expression of $\alpha_n, \beta_n,$ we have that
$$
X_n=(A_0\cdots A_{n-1})^{-1}\tau_{n+1},\quad S_n=\tau_n^{-1}A_0\cdots A_{n-1}, \quad n\geq0,\quad Y_0=\alpha_0\tau_0,
$$
and 
\begin{align*}
R_1&=\tau_1^{-1}(B_0-\alpha_0),& Y_1&=C_1(B_0-\alpha_0)^{-1}\tau_1,\\
R_2&=\tau_2^{-1}A_0(B_1-C_1(B_0-\alpha_0)^{-1}A_0),&Y_2&=C_2[B_1-C_1(B_0-\alpha_0)^{-1}A_0]^{-1}A_0^{-1}\tau_2,\quad\mbox{etc.}
\end{align*}

Therefore, if we are able to propose a good candidate for $\tau_n$, then we can compute all block entries $X_n,Y_n,R_n, S_n$ in terms of only $Y_0=\alpha_0\tau_0$.

\begin{remark}\label{Rems}
	In general, there are few conclusions we can derive for the sequence $\tau_n$ in terms of the positivity of the block entries $X_n,Y_n,R_n, S_n$. In particular, since $S_0=\tau_0^{-1}$, then $\tau_0^{-1}$, (but not necessarily $\tau_0$) must be a stochastic matrix. Also from $X_0=\tau_1$ and $Y_0=\alpha_0\tau_0$ we must have that $\tau_1$ and $\alpha_0\tau_0$ are semi-stochastic. In particular, since $\alpha_0=Y_0\tau_0^{-1}$ and $\tau_0^{-1}$ is stochastic, then $\alpha_0$ must be at least semi-stochastic as well.
\end{remark}

Similar considerations apply if we consider the LU block factorization \eqref{PxyLU}. Indeed, we will have
\begin{equation*}
P=\left[L\widetilde{\bm\beta}  \widetilde T\right]\left[\widetilde T^{-1}\widetilde{\bm\alpha}L^{-1}\right]=\widetilde P_L\widetilde P_U,
\end{equation*}
where now the coefficients $\widetilde\alpha_n,\widetilde\beta_n$ can be computed by
\begin{align*}
\widehat{B}_{n}&=\widetilde\beta_n+\widetilde\alpha_n,\quad n\geq0,\\
\widehat{C}_{n}&=\widetilde\beta_n\widetilde\alpha_{n-1},\quad n\geq1.
\end{align*}
Identifying block entries with \eqref{PxyLU} we get
\begin{align}
\label{xyrsgLU}\widetilde X_n&=\widetilde \tau_{n}^{-1}L_{n+1}^{-1},\quad \widetilde Y_n=\widetilde \tau_{n}^{-1}\widetilde\alpha_nL_n^{-1},\quad n\geq0,\\
\nonumber\widetilde S_n&=L_n\widetilde \tau_n,\quad \widetilde R_{n+1}=L_{n+1}\widetilde\beta_{n+1}\widetilde \tau_n,\quad n\geq0.
\end{align}
In this case we will get the relations
\begin{align*}
L_{n}(\widetilde\beta_{n}\widetilde \tau_{n-1}+\widetilde \tau_{n})\bm e_d&=\bm e_d,\quad n\geq0,\\
(\widetilde\alpha_{n}L_{n}^{-1}+L_{n+1}^{-1})\bm e_d&=\widetilde \tau_{n}\bm e_d,\quad n\geq0.
\end{align*}
Observe now that $\widetilde\tau_0$  must be stochastic ($\widetilde\beta_0$=0). This factorization will only depend on the sequence $\widetilde\tau_n$, since the sequences $\widetilde\alpha_n$ and $\widetilde\beta_n$ are uniquely determined in this case.

\section{Stochastic block Darboux transformations}\label{SEC3}

Assume that using the strategy above we have found appropriately $\alpha_0$ and $\tau_n$ such that $P_U$ and $P_L$ are stochastic matrices. We can perform what is called a \emph{discrete Darboux transformation} by reversing the order of multiplication. The Darboux transformation for second-order differential operators has a long history but probably the first reference of a discrete Darboux transformation like we study here appeared in \cite{MS} in connection with the Toda lattice.

If $P=P_UP_L$ as in \eqref{Pxy}, then by reversing the order of multiplication we obtain another block tridiagonal matrix of the form
\begin{equation}\label{DarbT}
\widetilde{P}=P_LP_U=\begin{pmatrix}
S_0&0&\\
R_1&S_1&0&\\
&\ddots&\ddots&\ddots
\end{pmatrix}\begin{pmatrix}
Y_0&X_0&\\
0&Y_1&X_1&\\
&\ddots&\ddots&\ddots
\end{pmatrix}=\begin{pmatrix}
\widetilde{B}_0&\widetilde{A}_0&\\
\widetilde{C}_1&\widetilde{B}_1&\widetilde{A}_1&\\
&\ddots&\ddots&\ddots
\end{pmatrix}.
\end{equation}
Now the new block entries are given by (see \eqref{alfbet} and \eqref{xyrsg})
\begin{align}
\nonumber\widetilde{A}_n&=S_nX_{n}=\tau_n^{-1}\tau_{n+1},\quad n\geq0,\\
\label{CoDT}\widetilde{B}_n&=R_nX_{n-1}+S_nY_n=\tau_n^{-1}(\beta_n+\alpha_n)\tau_{n}=\tau_n^{-1}(\widehat{B}_n-\beta_{n+1}+\beta_n)\tau_{n},\quad n\geq0,\\
\nonumber\widetilde{C}_n&=R_nY_{n-1}=\tau_n^{-1}(\beta_n\alpha_{n-1})\tau_{n-1}=\tau_n^{-1}(\beta_n\widehat{C}_{n-1}\beta_{n-1}^{-1})\tau_{n-1},\quad n\geq1.
\end{align}
The matrix $\widetilde{P}$ is actually stochastic, since the multiplication of two stochastic matrices is again a stochastic matrix. Therefore it is a new QBD process with block entries $\widetilde{A}_n$, $\widetilde{B}_n$ and $\widetilde{C}_n$. In fact we will get several QBD processes depending on many free parameters. In terms of a model driven by urn experiments (as we will see in Section \ref{SEC5}) the factorization $P=P_UP_L$ may be thought as two urn experiments, Experiment 1 and Experiment 2, respectively. We first perform the Experiment 1 and with the result we immediately perform the Experiment 2. The urn model for $\widetilde{P}=P_LP_U$ will proceed in the reversed order, first the Experiment 2 and with the result the Experiment 1.

The same can be done for the LU factorization \eqref{PxyLU} of the form $P=\widetilde P_L\widetilde P_U$. The corresponding Darboux transformation is
\begin{equation}\label{DarbT2}
\widehat{P}=\widetilde P_U\widetilde P_L=\begin{pmatrix}
\widetilde Y_0&\widetilde X_0&\\
0&\widetilde Y_1&\widetilde X_1&\\
&\ddots&\ddots&\ddots
\end{pmatrix}\begin{pmatrix}
\widetilde S_0&0&\\
\widetilde R_1&\widetilde S_1&0&\\
&\ddots&\ddots&\ddots
\end{pmatrix}=\begin{pmatrix}
\widehat{B}_0&\widehat{A}_0&\\
\widehat{C}_1&\widehat{B}_1&\widehat{A}_1&\\
&\ddots&\ddots&\ddots
\end{pmatrix}.
\end{equation}
The new coefficients are given by (see \eqref{xyrsgLU})
\begin{align*}
\widehat{A}_n&=\widetilde X_{n}\widetilde S_{n+1}=\widetilde\tau_n^{-1}\widetilde\tau_{n+1},\quad n\geq0,\\
\widehat{B}_n&=\widetilde X_{n}\widetilde R_{n+1}+\widetilde Y_n\widetilde S_n=\widetilde\tau_n^{-1}(\widetilde\beta_{n+1}+\widetilde\alpha_n)\widetilde\tau_n,\quad n\geq0,\\
\widehat{C}_n&=\widetilde Y_{n}\widetilde R_n=\widetilde\tau_n^{-1}\widetilde\alpha_n\widetilde\beta_{n}\widetilde\tau_{n-1},\quad n\geq1.
\end{align*}

If we assume that the block tridiagonal stochastic matrix $P$ is self-adjoint (in some appropriate Hilbert space) then there exists a unique weight matrix $W$ defined on the interval $-1\leq x\leq1$ (see \cite{DRSZ, G2}). Given such a weight matrix $W$ we can consider the skew symmetric bilinear form defined for any pair of matrix-valued polynomials $P(x)$ and $Q(x)$ by the numerical matrix
\begin{equation}\label{innp}
(P,Q)_{W}:=\int_{-1}^1P(x)W(x)Q^*(x)dx,
\end{equation}
where $Q^*(x)$ denotes the conjugate transpose of $Q(x)$. This leads, using the Gram-Schmidt process, to the existence of a sequence of matrix-valued orthogonal polynomials $Q_n$ with nonsingular leading coefficients satisfying  a three-term recursion relation
\begin{equation}\label{ttrrsts}
xQ_n(x)=A_{n}Q_{n+1}(x)+B_nQ_n(x)+C_{n}Q_{n-1}(x),\quad n\geq0,
\end{equation}
where $Q_{-1}(x)=0, Q_0(x)=I$ and $A_n,B_n,C_n$ are the coefficients of the block tridiagonal matrix $P$. With all these ingredients we can get the corresponding Karlin-McGregor representation formula for the block entry $(i,j)$ of $P^n$ (see \cite{DRSZ, G2}). Indeed,
\begin{equation}\label{KMcGRF}
P_{ij}^n=(x^nQ_i,Q_j)_W(Q_j,Q_j)_W^{-1}=\left(\int_{-1}^1x^nQ_i(x)W(x)Q_j^*(x)dx\right)\left(\int_{-1}^1Q_j(x)W(x)Q_j^*(x)dx\right)^{-1}.
\end{equation}
We can also compute the invariant measure $\bm\pi$ of the QBD process $P$ in terms of the inverse of the norms $\Pi_n=(Q_n,Q_n)_W^{-1}$ (see Theorem 3.1 of \cite{dI2}). Indeed, this invariant vector is given by
\begin{equation}\label{InvMea}
\bm\pi=\left((\Pi_0\bm e_d)^T; (\Pi_1\bm e_d)^T; (\Pi_2\bm e_d)^T;\cdots\right),
\end{equation}
where we recall here that $\bm e_d$ is the $d$-dimensional vector with all components equal to 1.

One important aspect of the Darboux transformation starting from the UL factorization is to study how to transform the \emph{matrix-valued} spectral measure associated with a QBD process with one-step transition probability matrix $P$. The property of being self-adjoint may be lost for the Darboux transformation of $P$ given by $\widetilde P$.

In the scalar case (tridiagonal matrix $P$) it is very well known that if the moment $\mu_{-1}=\int_{-1}^1d\omega(x)/x$ is well defined, where $\omega$ is the spectral measure associated with $P$, then a candidate for the family of spectral measures is then given by a \emph{Geronimus transformation} of $\omega$, i.e.
\begin{equation*}\label{spmes}
\widetilde\omega(x)=y_0\frac{\omega(x)}{x}+M\delta_0(x),\quad M=1-y_0\mu_{-1},
\end{equation*}
where $\delta_0(x)$ is the Dirac delta located at $x=0$. Similarly, for the LU factorization, the corresponding Darboux transformation \eqref{DarbT2} $\widehat{P}$ gives rise a spectral measure $\widehat\omega$ given by a \emph{Christoffel transformation} of $\omega$, i.e. $\widehat\omega(x)=x\omega(x)$ (see \cite{GdI3} and references therein).

In the matrix-valued case it is possible to see (in analogy with the scalar case), that if the moment $\mu_{-1}=\int_{-1}^1dW(x)/x$ is well defined, then a candidate for the family of matrix-valued spectral measures associated with the Darboux transformation $\widetilde P$ is again a Geronimus transformation of $W$, i.e.
\begin{equation}\label{spme}
\widetilde W(x)=\frac{W(x)}{x}+\bm M\delta_0(x),\quad \bm M=\alpha_0^{-1}\mu_0-\mu_{-1},
\end{equation}
as we will see in the example of the next section. Observe from the derivation of the coefficients $\widetilde{A}_n,\widetilde{B}_n,\widetilde{C}_n$ in \eqref{CoDT}, that the free parameters of $\widetilde W$ will only depend on $\alpha_0$ and not on the sequence $\tau_n$, which will only interfere in the normalization of the corresponding matrix-valued polynomials. In the case where $\alpha_0$ is a singular matrix, we will have a degenerate matrix-valued spectral measure. Also we observe that $\widetilde W$ in general \emph{is neither symmetric nor positive semidefinite}. 
In order for $\widetilde W$ to be a proper weight matrix, the matrix $\bm M$ in \eqref{spme} has to be positive semidefinite.

Similarly, for the LU factorization, the corresponding Darboux transformation \eqref{DarbT2} $\widehat{P}$ gives rise to a block Jacobi matrix and a matrix-valued spectral measure $\widehat W$. It is possible to see that $\widehat W$ is given by a Christoffel transformation of $W$, i.e.
\begin{equation*}\label{spme2}
\widehat W(x)=x W(x).
\end{equation*}
In this case the weight matrix $\widehat W$ is unique and positive semidefinite.

\section{The (Jacobi type) one-step example}\label{SEC4}

In this section we will study a specific example coming from group representation theory, introduced for the first time in \cite{GPT1}. In \cite{GdI2} we studied the probabilistic aspects of this example and gave an explicit expression of the block entries $A_n, B_n, C_n$ in \eqref{PP}. The most general situation is considered  in \cite{GPT3}, where the authors also give two stochastic models in terms of urns and Young diagrams.

First we will rewrite the block entries $A_n, B_n, C_n$ introduced in \cite{GdI2} in a more convenient way (following the ideas of \cite{GPT3}, see also \cite{GPT2}). For $\alpha,\beta>-1$ and $0<k<\beta+1$ define the coefficients
\begin{align*}
a_1(i,n)&=\frac{(n+k)(n+\beta+d)}{(2n+\alpha+\beta+d+i)(n+k+d-i-1)},\hspace{.95cm} b_1(i,n)=\frac{n(n+k+d-1)}{(2n+\alpha+\beta+d+i-1)(n+k+d-i-1)},\\
a_2(i,n)&=\frac{(d-i-1)(\beta-k+i+1)}{(n+\alpha+\beta-k+2i+1)(n+k+d-i-1)},\quad b_2(i,n)=\frac{i(k+d-i-1)}{(n+\alpha+\beta-k+2i)(n+k+d-i-1)},\\
a_3(i,n)&=\frac{(n+\alpha+i)(n+\alpha+\beta-k+d+i)}{(2n+\alpha+\beta+d+i)(n+\alpha+\beta-k+2i+1)},\; b_3(i,n)=\frac{(n+\alpha+\beta+d+i-1)(n+\alpha+\beta-k+i)}{(2n+\alpha+\beta+d+i-1)(n+\alpha+\beta-k+2i)}.
\end{align*}
Observe that all these coefficients are always \emph{positive} for $i=0,1,\ldots, d-1$, $n\geq0$, and that we have
\begin{equation}\label{asbs}
a_1(i,n)+a_2(i,n)+a_3(i,n)=1,\quad b_1(i,n)+b_2(i,n)+b_3(i,n)=1.
\end{equation}
Observe also that when $d\geq2$ there is a new parameter $k$. The block entries $A_n, B_n, C_n$ are given by
\begin{align*}
A_n&=\sum_{i=0}^{d-1}a_1(i,n)b_3(i,n+1)E_{ii}+\sum_{i=0}^{d-2}a_1(i+1,n)b_2(i+1,n+1)E_{i+1,i},\\
B_n&=\sum_{i=0}^{d-2}a_3(i+1,n)b_2(i+1,n)E_{i+1,i}+\sum_{i=0}^{d-2}a_2(i,n)b_3(i+1,n)E_{i,i+1}\\
&\quad+\sum_{i=0}^{d-1}\left[a_1(i,n)b_1(i,n+1)+a_2(i,n)b_2(i+1,n)+a_3(i,n)b_3(i,n)\right]E_{ii},\\
C_n&=\sum_{i=0}^{d-1}a_3(i,n)b_1(i,n)E_{ii}+\sum_{i=0}^{d-2}a_2(i,n)b_1(i+1,n)E_{i,i+1},
\end{align*}
where $E_{ij}$ denotes as usual the matrix with entry $(i,j)$, which is equal to 1 and 0 elsewhere.

The family of matrix-valued polynomials generated by the three-term recursion relation \eqref{ttrrsts}, where $Q_{-1}(x)=0$ and $Q_0(x)=I$, is in fact orthogonal with respect to the weight matrix (see \cite{PT1} and \cite{GdI2})
\begin{equation}\label{WW}
W(x)=x^{\alpha}(1-x)^{\beta}V^*Z(x)V,\quad x\in[0,1],
\end{equation}
where
\begin{equation*}\label{ZZ}
Z(x)=\sum_{i,j=0}^{d-1}\bigg(\sum_{r=0}^{d-1}\begin{pmatrix}
                                       r \\
                                       i \\
                                     \end{pmatrix}\begin{pmatrix}
                                       r \\
                                       j \\
                                     \end{pmatrix}\begin{pmatrix}
                                       d+k-r-2 \\
                                       d-r-1 \\
                                     \end{pmatrix}\begin{pmatrix}
                                      \beta-k+r \\
                                       r \\
                                     \end{pmatrix}(1-x)^{i+j}x^{d-r-1}\bigg)E_{ij},
\end{equation*}
and $V$ is the nonsingular upper triangular matrix
\begin{equation*}\label{P}
V=\sum_{i\leq j}
(-1)^{i}\frac{(-j)_{i}}{(1-d)_{i}}\frac{(\alpha+\beta-k+j+1)_{i}}{(\beta-k+1)_{i}}E_{ij}.
\end{equation*}
Here $(a)_n=a(a+1)\cdots (a+n-1)$ denotes the Pochhammer symbol.

The way of writing the block entries $A_n, B_n, C_n$ above is very useful when trying to find a particular factorization of the form $P=P_UP_L$ as in \eqref{Pxy}. Indeed, good candidates for the block entries $X_n, Y_n, R_n, S_n$ are given by
\begin{align}
\label{xyrs1}X_n&=\sum_{i=0}^{d-1}a_1(i,n)E_{ii},\quad Y_n=\sum_{i=0}^{d-1}a_3(i,n)E_{ii}+\sum_{i=0}^{d-2}a_2(i,n)E_{i,i+1},\\
\label{xyrs2}R_n&=\sum_{i=0}^{d-1}b_1(i,n)E_{ii},\quad S_n=\sum_{i=0}^{d-1}b_3(i,n)E_{ii}+\sum_{i=0}^{d-2}b_2(i+1,n)E_{i+1,i}.
\end{align}

The relations \eqref{asbs} actually say that $P_U$ and $P_L$ are stochastic matrices. This factorization appeared for the first time in \cite{GPT3}, but this \emph{is not the only} UL stochastic factorization possible for the matrix $P$. As we saw in the previous section, the UL factorization comes with at least a free extra parameter $\alpha_0$ which in this case is a $d\times d$ matrix with some restrictions.

We will see now how the choice of factors in \cite{GPT3} fits with the more general framework above, i.e. we
give now a choice for $\alpha_0$ and $\tau_n$ such that we get $X_n, Y_n, R_n, S_n$ as in \eqref{xyrs1} and \eqref{xyrs2}, i.e. as in \cite{GPT3}.

Consider
\begin{equation}\label{alf0}
\alpha_0=B_0-D_0,
\end{equation}
where $D_0$ is the diagonal matrix
$$
D_0=\sum_{i=0}^{d-1}a_1(i,0)b_1(i,1)E_{ii}=\sum_{i=0}^{d-1}\frac{k(\beta+d)(k+d)}{(\alpha+\beta+d+i)_2(k+d-i-1)_2}E_{ii}.
$$
Observe that $D_0$ is the first summand in the diagonal entries of $B_0$. Let $L_n$ be as in \eqref{LLn}\footnote{A different way of writing the lower triangular matrix $L_n$ can be found in page 751 of \cite{GdI2} (written there as $A_n^n$). } and choose
\begin{equation*}\label{taun}
\tau_n=\tau_0\left(\left.L_n^{-1}\right|_{\alpha=\alpha-1}\right),
\end{equation*}
where $\tau_0$ is a lower triangular matrix which inverse is given by the stochastic matrix
$$
\tau_0^{-1}=\sum_{i=0}^{d-1}\frac{\alpha+\beta-k+i}{\alpha+\beta-k+2i}E_{ii}+\sum_{i=0}^{d-2}\frac{i+1}{\alpha+\beta-k+2i+2}E_{i+1,i}=S_0.
$$
Then the block entries $X_n, Y_n, R_n, S_n$ in \eqref{xyrsg} coincide with the block entries $X_n, Y_n, R_n, S_n$ given in \cite{GPT3} (see \eqref{xyrs1} and \eqref{xyrs2}). Moreover, the matrix-valued spectral measure $\widetilde W$ associated with the Darboux transformation $\widetilde P=P_LP_U$ (see \eqref{DarbT}) is given by
$$
\widetilde W(x)=\frac{W(x)}{x},\quad x\in[0,1],\quad \alpha>0,\quad \beta>-1,\quad 0<k<\beta+1,
$$
where $W$ is the matrix-valued spectral measure associated with $P$ (see \eqref{WW}). Observe that $\alpha_0$ in this case is chosen in such a way that $ \bm M=\alpha_0^{-1}\mu_0-\mu_{-1}=0$ in \eqref{spme}, i.e. $\alpha_0=\mu_0\mu_{-1}^{-1}$, an alternative to the expression for $\alpha_0$ above.

\medskip

If we consider the LU factorization of the same block tridiagonal $P$, let us choose
\begin{equation*}\label{taun2}
\widetilde\tau_n=\widetilde\tau_0\left(\left.L_n^{-1}\right|_{\alpha=\alpha+1}\right),
\end{equation*}
where $\widetilde\tau_0=\left.\tau_0^{-1}\right|_{\alpha=\alpha+1}$. Then the block entries $\widetilde X_n, \widetilde Y_n, \widetilde R_n, \widetilde S_n$ in \eqref{xyrsgLU} are given by
\begin{equation*}
\widetilde X_n=\left.X_n\right|_{\alpha=\alpha+1},\quad \widetilde Y_n=\left.Y_n\right|_{\alpha=\alpha+1},\quad \widetilde R_n=\left.R_n\right|_{\alpha=\alpha+1},\quad \widetilde S_n=\left.S_n\right|_{\alpha=\alpha+1}.
\end{equation*}
The matrix-valued spectral measure $\widehat W$ associated with the Darboux transformation $\widehat P=\widetilde P_U\widetilde P_L$ of $P$ (see \eqref{DarbT2}) is given by the Christoffel transformation of $W$, i.e.
\begin{equation*}
\widehat W(x)=x W(x),\quad x\in[0,1],\quad \alpha,\beta>-1,\quad 0<k<\beta+1.
\end{equation*}

\medskip

As we said earlier, the choice of $\alpha_0$ above is the one that among the possible factorizations of $P$ reproduces the results in \cite{GPT3}.

\bigskip

Let us focus now in the case where $\alpha_0$ is not necessarily chosen as in \eqref{alf0}. For simplicity, we will explore only the case where $d=2$. The case $d=1$ was already studied in \cite{GdI3}, along with an urn model for the Jacobi polynomials (for a different urn model for the Jacobi polynomials which does not exploit the factorization of the tridiagonal matrix see \cite{G4}).

\subsection{Case $d=2$}

The coefficients of $P$ given earlier become now
\begin{align}
\label{aan}A_n&=\begin{pmatrix}
\frac{(\beta+n+2) (k+n) (\alpha+\beta+n+2)}{(k+n+1) (\alpha+\beta+2 n+2) (\alpha+\beta+2 n+3)}& 0\\
\frac{k (\beta+n+2)}{(\alpha+\beta-k+n+3) (\alpha+\beta+2 n+3) (k+n+1)}&\frac{(\beta+n+2) (\alpha+\beta+n+3) (\alpha+\beta-k+n+2)}{(\alpha+\beta+2 n+3) (\alpha+\beta+2 n+4) (\alpha+\beta-k+n+3)}
\end{pmatrix},\\
\label{bbn}B_n&=\begin{pmatrix}
B_n^{11}& \frac{(\beta-k+1)(\alpha+\beta+n+2)}{(k+n+1)(\alpha+\beta+2n+2)(\alpha+\beta-k+n+2)}\\
\frac{(\alpha+n+1)k}{(k+n)(\alpha+\beta-k+n+2)(\alpha+\beta+2n+3)}&B_n^{22}
\end{pmatrix},\\
\label{ccn}C_n&=\begin{pmatrix}
\frac{n (\alpha+n) (\alpha+\beta-k+n+2)}{(\alpha+\beta-k+n+1) (\alpha+\beta+2 n+1) (\alpha+\beta+2 n+2)}& \frac{n (\beta-k+1)}{(\alpha+\beta-k+n+1) (\alpha+\beta+2 n+2) (k+n)}\\
0&\frac{n (\alpha+n+1) (k+n+1)}{(k+n) (\alpha+\beta+2 n+2) (\alpha+\beta+2 n+3)}
\end{pmatrix},
\end{align}
where
\begin{align*}
B_n^{11}&=\frac{(n+k)(n+\beta+2)(n+1)}{(\alpha+\beta+2n+2)(n+k+1)(\alpha+\beta+2n+3)}+\frac{(n+\alpha)(\alpha+\beta-k+n+2)(n+\alpha+\beta+1)}{(\alpha+\beta+2n+2)(\alpha+1+n-k+\beta)(\alpha+\beta+2n+1)}\\
&\quad+\frac{k(\beta-k+1)}{(\alpha+1+n-k+\beta)(n+k+1)(\alpha+\beta-k+n+2)(n+k)},
\end{align*}
and
\begin{align*}
B_n^{22}&=\frac{(n+\beta+2)(n+1)(n+k+2)}{(\alpha+\beta+2n+3)(\alpha+\beta+2n+4)(n+k+1)}+\frac{(\alpha+n+1)(\alpha+\beta+n+2)(\alpha+1+n-k+\beta)}{(\alpha+\beta+2n+3)(\alpha+\beta+2n+2)(\alpha+\beta-k+n+2)}.
\end{align*}
The main difference with the scalar case ($d=1$) is that we have a new parameter $0<k<\beta+1$. A straightforward computation, using the definition of $a_j(i,n),b_j(i,n), j=1,2,3,$ gives
$$
L_n=\begin{pmatrix}
\frac{(n+k)(\alpha+\beta+n+2)_n}{k(\beta+2)_n}&0\\
-\frac{n(\alpha+\beta+n+3)_n}{(\alpha+\beta-k+2)(\beta+2)_n}&\frac{(\alpha+\beta+n-k+2)(\alpha+\beta+n+3)_n}{(\alpha+\beta-k+2)(\beta+2)_n}
\end{pmatrix}.
$$
The coefficients $\alpha_n,\beta_n$ of the UL factorization \eqref{Jab} are going to be computed by using the expressions \eqref{alfbet}
Since the matrix $\alpha_0$ is in principle any $2\times2$ matrix, the coefficients $\alpha_n,\beta_n$ are not easy to find. Once we have these coefficients we can compute $\tau_n$ by solving \eqref{eq2} and \eqref{eq3}, which now also may yield more degrees of freedom.

\medskip

We first review in this $d=2$ case what we did earlier for general $d$.

\medskip

If we choose $\alpha_0$ as in \eqref{alf0}, i.e.
\begin{align}
\label{alf0m}\alpha_0=\begin{pmatrix}
\frac{\beta-k+1}{(1+\alpha+\beta-k)(1+k)(2+\alpha+\beta-k)}+\frac{\alpha(2+\alpha+\beta-k)}{(2+\alpha+\beta)(1+\alpha+\beta-k)}&\frac{\beta-k+1}{(1+k)(2+\alpha+\beta-k)}\\
\frac{1+\alpha}{(3+\alpha+\beta)(2+\alpha+\beta-k)}&\frac{(1+\alpha)(1+\alpha+\beta-k)}{(3+\alpha+\beta)(2+\alpha+\beta-k)}
\end{pmatrix},
\end{align}
then the explicit expression for $\tau_n$ is given by
$$
\tau_0^{-1}\tau_n=\begin{pmatrix}
\frac{k(\beta+2)_n}{(n+k)(\alpha+\beta+n+1)_n}&0\\
\frac{n(\beta+2)_n}{(n+k)(\alpha+\beta+n-k+1)(\alpha+\beta+n+1)_n}&\frac{(\alpha+\beta-k+1)(\beta+2)_n}{(\alpha+\beta+n-k+1)(\alpha+\beta+n+2)_n}
\end{pmatrix}=\left(\left.L_n^{-1}\right|_{\alpha=\alpha-1}\right),
$$
where
\begin{equation}\label{tau0}
\tau_0^{-1}=\begin{pmatrix}
1&0\\
\frac{1}{\alpha+\beta-k+2}&\frac{\alpha+\beta-k+1}{\alpha+\beta-k+2}
\end{pmatrix}.
\end{equation}
The block entries of the stochastic matrices $P_U$ and $P_L$ are given by
\begin{align}
\label{Xn2}
X_n&=\begin{pmatrix}
\frac{(n+k)(n+\beta+2)}{(2n+\alpha+\beta+2)(n+k+1)}&0\\
0&\frac{n+\beta+2}{2n+\alpha+\beta+3}
\end{pmatrix},\quad 
Y_n=\begin{pmatrix}
\frac{(n+\alpha)(n+\alpha+\beta-k+2)}{(2n+\alpha+\beta+2)(n+\alpha+1-k+\beta)}&\frac{\beta-k+1}{(n+\alpha+1-k+\beta)(n+k+1)}\\
0&\frac{n+\alpha+1}{2n+\alpha+\beta+3}
\end{pmatrix},\\
\label{Sn2}
S_n&=\begin{pmatrix}
\frac{n+\alpha+\beta+1}{2n+\alpha+\beta+1}&0\\
\frac{k}{(n+\alpha+\beta-k+2)(n+k)}&\frac{(n+\alpha+\beta+2)(n+\alpha+1-k+\beta)}{(2n+\alpha+\beta+2)(n+\alpha+\beta-k+2)}
\end{pmatrix},\quad 
R_n=\begin{pmatrix}
\frac{n}{2n+\alpha+\beta+1}&0\\
0&\frac{n(n+k+1)}{(2n+\alpha+\beta+2)(n+k)}
\end{pmatrix}.
\end{align}
We will give in the next section an interpretation of these matrices in terms of an urn model.

\bigskip

We are done with trying to reproduce the resuls in \cite{GPT3} and we move on to a generic $\alpha_0$ where things are more complicated.
The only thing that we know is that $\alpha_0$ must be at least semi-stochastic (see Remark \ref{Rems}).

\medskip

Getting away from \eqref{alf0m}, let us put
$$
\alpha_0=\begin{pmatrix}s_{11}&s_{12}\\s_{21}&s_{22}\end{pmatrix}.
$$
According to \eqref{spme}, the family of matrix-valued spectral measures associated with the Darboux transformation $\widetilde P$ is given by
\begin{equation}\label{WWg2}
\widetilde W(x)=\frac{W(x)}{x}+\bm M\delta_0(x),\quad \bm M=\alpha_0^{-1}\mu_0-\mu_{-1}.
\end{equation}
In the case of this example we have that
$$
\mu_0=\frac{\Gamma(\alpha+1)\Gamma(\beta+2)(\alpha+\beta-k+2)}{\Gamma(\alpha+\beta+3)}\begin{pmatrix} 1 &0 \\
0&\frac{(\alpha+1)(k+1)}{(\alpha+\beta+3)(\beta-k+1)} \end{pmatrix},
$$
and (assuming $\alpha>0, \beta>-1$)
$$
\mu_{-1}=\frac{\Gamma(\alpha)\Gamma(\beta+2)}{\Gamma(\alpha+\beta+2)}\begin{pmatrix} \alpha+\beta-k+1 &-1 \\
-1&\frac{(\alpha+1)(k+1)(\alpha+\beta-k+2)-k(\beta-k+1)}{(\alpha+\beta+2)(\beta-k+1)} \end{pmatrix}.
$$
Observe that we have the following relation between the moments $\mu_0$ and $\mu_{-1}$:
$$
\left.\mu_0\right|_{\alpha=\alpha-1}=\tau_0^{-1}\mu_{-1}\tau_0^{-*},
$$
where $\tau_0^{-1}$ is the lower triangular matrix \eqref{tau0}. As we mentioned in Section \ref{SEC3}, $\widetilde W$ is in general neither symmetric nor positive semidefinite. It is easy to see that $\bm M$ in \eqref{WWg2} is symmetric if and only if one of the entries of $\alpha_0$ is chosen according to the following relation
\begin{equation}\label{s12s}
s_{12}=\frac{(\beta-k+1)(\alpha+\beta+3)}{(\alpha+1)(k+1)}s_{21}. 
\end{equation}

In what follows we will study two special cases where we can explicitly analyze the different values of the parameters $s_{ij}$ for which the Darboux transformation gives rise to a QBD process. In the first case we will focus on the positivity of the matrix-valued spectral measure $\widetilde W$ in \eqref{WWg2}, while in the second case we will analyze the stochastic block factorization without looking at the positivity or the symmetry of $\widetilde W$.

\begin{enumerate}

\item Let us choose for convenience
\begin{equation}\label{alf0s}
\alpha_0=\begin{pmatrix}\frac{(\alpha+\beta+3)(\beta-k+1)}{(k+1)(\alpha+1)(\alpha+\beta-k+1)}s_{11}&\frac{(\alpha+\beta+3)(\beta-k+1)}{(k+1)(\alpha+1)}s_{21}\\s_{21}&(\alpha+\beta-k+1)s_{21}\end{pmatrix},
\end{equation}
		where $s_{11}$ and $s_{21}$ are in principle free parameters. Since condition \eqref{s12s} is satisfied, $\widetilde W$ in \eqref{WWg2} is a symmetric matrix. The determinant of $\alpha_0$ is given by
$$
|\alpha_0|=\frac{(\alpha+\beta+3)(\beta-k+1)}{(k+1)(\alpha+1)}s_{21}(s_{11}-s_{21}).
$$
A straightforward computation shows that $\bm M$ can be written in this case as
$$
\bm M=\tau_0\begin{pmatrix}\frac{(\alpha)_2(k+1)(\alpha+\beta-k+2)}{(\beta-k+1)(\alpha+\beta+2)_2(s_{11}-s_{21})}-1&0 \\
0&\frac{\alpha+1}{(\alpha+\beta+2)(\alpha+\beta-k+2)s_{21}}-1\end{pmatrix}\left.\mu_0\right|_{\alpha=\alpha-1}\tau_0^*,
$$
where $\tau_0$ is the inverse of the the lower triangular matrix \eqref{tau0}. Since $\mu_0$ is diagonal, then $\bm M$ is positive semidefinite if the parameters $s_{11}$ and $s_{21}$ are chosen in the following range
\begin{align}
\label{upboun}0<&s_{21}\leq\frac{\alpha+1}{(\alpha+\beta+3)(\alpha+\beta-k+2)},\\
\nonumber s_{21}<&s_{11}\leq s_{21}+\frac{(\alpha)_2(k+1)(\alpha+\beta-k+2)}{(\beta-k+1)(\alpha+\beta+2)_2}.
\end{align}

\medskip

Based on what happens in the scalar case, one could expect that if we choose $s_{11}$ and $s_{21}$ in the range above, then we should get a stochastic block factorization of $P$ of the form \eqref{Pxy}, where $P_U$ and $P_L$ are also stochastic matrices. However, \emph{this is not true}. In fact the condition on the positivity of the entries of $X_n,Y_n,S_n, R_n$ will force us to modify the upper bound of the second inequality in \eqref{upboun}.

Another important point is the choice of the sequence of matrices $\tau_n$ such that \eqref{eq2} and \eqref{eq3} hold and the entries of the matrices $X_n,Y_n,S_n, R_n$ in \eqref{xyrsg} are all nonnegative. Since $\tau_n$ is a $2\times2$ matrix, it has 4 degrees of freedom for every $n$. For simplicity, we look for \emph{lower triangular matrices} $\tau_n$ with 3 degrees of freedom for every $n$. The diagonal entries of $\tau_n$ can be given by solving the equations \eqref{eq2} and \eqref{eq3}. Then, the matrices $X_n$ and $S_n$ are always lower triangular, while the matrices $Y_n$ and $R_n$ are both full matrices. In order to surmise the remaining free parameter of $\tau_n$, we force $Y_n$ to be upper triangular and as a consequence $R_n$ will be also upper triangular. These conditions will give us explicitly $\tau_n$ for every $n$ and one sees that the sum of every row of $P_U$ and $P_L$ is always 1.

Finally, we need all entries of $X_n,Y_n,S_n, R_n$ to be nonnegative. The entries of these matrices are rational functions depending on $s_{11}$ and $s_{21}$. After extensive symbolic computations we find that all entries of $X_n,Y_n,S_n, R_n$ are nonnegative (and therefore $P_U$ and $P_L$ are stochastic matrices) if the parameters $s_{11}$ and $s_{21}$ are chosen in the following range
\begin{align*}
0<&s_{21}\leq\frac{\alpha+1}{(\alpha+\beta+3)(\alpha+\beta-k+2)},\\
s_{21}<&s_{11}\leq \frac{s_{21}\left(s_{21}-\frac{(\alpha+1)^2(k+1)}{k(\beta-k+1)(\alpha+\beta+3)}\right)}{s_{21}-\frac{(\alpha+1)(k+1)}{k(\alpha+\beta-k+1)(\alpha+\beta+3)}}.
\end{align*}
Observe that the singular point of the rational function of $s_{21}$ above is to the right of the upper bound of $s_{21}$ in \eqref{upboun}. In Figure \ref{fig1} we can see this region for the special values of $\alpha=3, \beta=2, k=1$. The green line (above the shaded area of the figure) is the upper bound for which $\bm M$ is positive semidefinite, but we observe here that there may be values of $s_{11},s_{21}$ for which $\bm M$ is positive semidefinite, and yet the block entries $X_n,Y_n,S_n, R_n$ do not have all their entries nonnegative.

\begin{figure}[h]
\begin{center}
\vspace{-0.0cm}
\includegraphics[height=7.4cm]{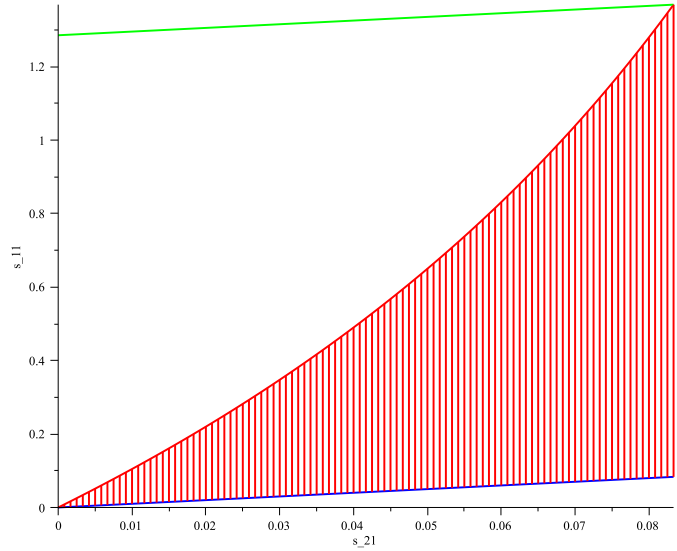}
\end{center}
\vspace{-0.2cm}
	\caption{The region with red stripes (shaded area) gives all possible values of $s_{21}$ and $s_{11}$ for which all entries of $X_n,Y_n,S_n, R_n$ are nonnegative for the values of $\alpha=3, \beta=2, k=1$. The green line (above the shaded area) is the upper bound for which $\bm M$ is positive semidefinite.}
\label{fig1}
\end{figure}

This concludes our look at the case when $\alpha_0$ was chosen as in \eqref{alf0s}.

\item Let us focus now exclusively on the block entries $X_n,Y_n,S_n, R_n$ and the sequence of matrices $\tau_n$ chosen to guarantee that $P_U$ and $P_L$ are stochastic without looking into the matrix-valued spectral measure $\widetilde W$ resulting after the Darboux transformation. It was mentioned in the previous case that we can always choose a unique sequence of lower triangular matrices $\tau_n$ such that \eqref{eq2} and \eqref{eq3} hold and $X_n$ and $S_n$ are lower triangular matrices while $Y_n$ and $R_n$ are upper triangular matrices. Imagine now that we would like to have one (or several) of the matrices $X_n,Y_n,S_n, R_n$ \emph{as a diagonal one}. To insure this we need to impose some restrictions on the parameters $s_{ij}$ of $\alpha_0$. There are four possible situations:
\begin{enumerate}
	\item $X_n$ \emph{diagonal}: the matrix $\alpha_0$ should be chosen as the two-parameter family of matrices
$$
\alpha_0=\begin{pmatrix}
s_{11}&s_{12}\\
\frac{1+\alpha}{(3+\alpha+\beta)(2+\alpha+\beta-k)}&\frac{(1+\alpha)(1+\alpha+\beta-k)}{(3+\alpha+\beta)(2+\alpha+\beta-k)}
\end{pmatrix}.
$$
Observe that the second row is the same as the second row of $\alpha_0$ in \eqref{alf0m}.
\item $R_n$ \emph{diagonal}: the matrix $\alpha_0$ should be chosen as the two-parameter family of matrices
$$
\alpha_0=\begin{pmatrix}
s_{11}&\frac{\beta-k+1}{(1+k)(2+\alpha+\beta-k)}\\
s_{21}&\frac{(1+\alpha)(1+\alpha+\beta-k)}{(3+\alpha+\beta)(2+\alpha+\beta-k)}
\end{pmatrix}.
$$
Observe that the second column is the same as the second column of $\alpha_0$ in \eqref{alf0m}.
\item $Y_n$ \emph{diagonal}: the matrix $\alpha_0$ should be chosen as the two-parameter family of singular matrices
$$
\alpha_0=\begin{pmatrix}
0&0\\
s_{21}&s_{22}
\end{pmatrix}.
$$
\item $S_n$ \emph{diagonal}: the matrix $\alpha_0$ should be chosen as the two-parameter family of singular matrices
$$
\alpha_0=\begin{pmatrix}
0&s_{21}\\
0&s_{22}
\end{pmatrix}.
$$
\end{enumerate}
Let us analyze further the case (a) (the rest can be studied in a similar manner). For convenience we will work with the normalization
$$
\alpha_0=\begin{pmatrix}
\left[\frac{\beta-k+1}{(1+\alpha+\beta-k)(1+k)(2+\alpha+\beta-k)}+\frac{\alpha(2+\alpha+\beta-k)}{(2+\alpha+\beta)(1+\alpha+\beta-k)}\right]s_{11}&\frac{\beta-k+1}{(1+k)(2+\alpha+\beta-k)}s_{12}\\
\frac{1+\alpha}{(3+\alpha+\beta)(2+\alpha+\beta-k)}&\frac{(1+\alpha)(1+\alpha+\beta-k)}{(3+\alpha+\beta)(2+\alpha+\beta-k)}
\end{pmatrix}.
$$
Then, all entries of $X_n,Y_n,S_n, R_n$ are nonnegative (and therefore $P_U$ and $P_L$ are stochastic matrices) if $s_{11}$ and $s_{12}$ are chosen in the following range
\begin{align*}
0<&s_{11}\leq1,\\
s_{11}\leq &s_{12}\leq \min\left\{\left(1+\frac{\alpha(k+1)(\alpha+\beta-k+2)^2}{(\alpha+\beta+2)(\beta-k+1)}\right)s_{11},1\right\}.
\end{align*}
In Figure \ref{fig2} we can visualize this region for the special values of $\alpha=1, \beta=2, k=2$.
\begin{figure}[h]
\begin{center}
\vspace{-0.0cm}
\includegraphics[height=7.4cm]{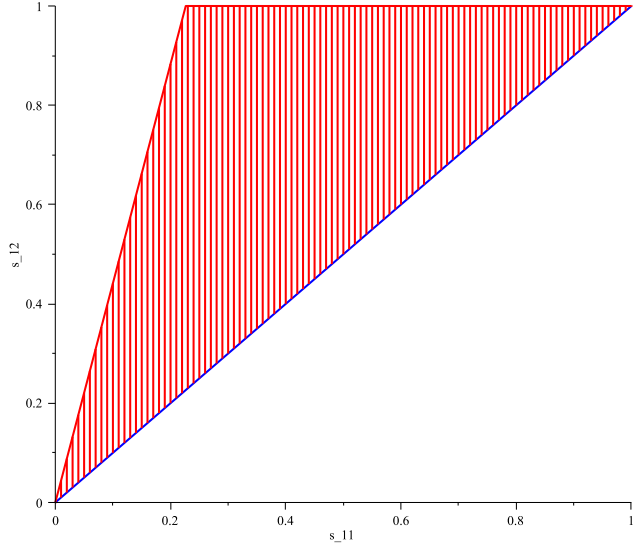}
\end{center}
\vspace{-0.2cm}
	\caption{The region with red stripes (shaded area) gives all possible values of $s_{11}$ and $s_{12}$ for which all entries of $X_n,Y_n,S_n, R_n$ are nonnegative for the values of $\alpha=1, \beta=2, k=2$.}
\label{fig2}
\end{figure}
Since in general the condition \eqref{s12s} is not satisfied, then there is no symmetric and positive definite spectral measure $\widetilde W$ of the form \eqref{WWg2} for the Darboux transformation $\widetilde P$. Nevertheless, if we consider the monic matrix-valued polynomials $\widetilde P_n$ generated by the Jacobi matrix $\widetilde J=\bm\beta\bm\alpha$ (see \eqref{Jab}), then they are \emph{left-orthogonal} (i.e. $(\widetilde P_n,\widetilde P_m)_{\widetilde W}=0$ for $n>m$, see \eqref{innp}) with respect to the matrix-valued function \eqref{WWg2}, with
$$
\bm M=\tau_0\begin{pmatrix}\frac{(\alpha)_2}{(\alpha+\beta+2)_2|\alpha_0|}-1&\frac{(\beta-k+1)(\alpha+1)(1-s_{12})}{|\alpha_0|(k+1)(\alpha+\beta+3)(\alpha+\beta-k+2)}\\
0&0\end{pmatrix}\left.\mu_0\right|_{\alpha=\alpha-1}\tau_0^*,
$$
where $\tau_0$ is the inverse of the the lower triangular matrix \eqref{tau0}. Observe that $\bm M$ is a singular and non-symmetric matrix. Observe also that the only case where $\bm M$ is symmetric is given by choosing $s_{12}=1$. But this is just the first case obtained by taking
$$
s_{21}=\frac{\alpha+1}{(\alpha+\beta+2)(\alpha+\beta-k+2)}.
$$
\begin{remark}
	We would like to remark on a special property of the matrix-valued orthogonal polynomials $\widetilde P_n$ generated by the Darboux transformation of $P$ for the case $s_{12}=1$ above. It is well known that the original matrix-valued orthogonal polynomials $P_n$ satisfy a second-order differential equation of the form
\begin{equation}\label{sode}
P_n''(x)F_2(x)+P_n'(x)F_1(x)+P_n(x)F_0=\Lambda_nP_n(x),
\end{equation}
where $F_2(x)=x(1-x)I$ and $F_1, F_0$ certain matrix polynomials of degree 1 and 0, respectively (see for instance \cite{PT1} or \cite{GdI2}). In this situation (and only in this situation) the matrix-valued polynomials $\widetilde P_n$ obtained by performing the Darboux transformation also satisfy a \emph{second-order differential equation} of the form \eqref{sode} with coefficients $\widetilde F_2, \widetilde F_1, \widetilde F_0$ given by
\begin{align*}
\widetilde F_2(x)&=x\left[\begin{pmatrix}0&0\\1&-1\end{pmatrix}x+\begin{pmatrix}\frac{\beta-k+1}{\alpha+\beta-k+2}&-\frac{\beta-k+1}{\alpha+\beta-k+2}\\-\frac{\alpha+1}{\alpha+\beta-k+2}&\frac{\alpha+1}{\alpha+\beta-k+2}\end{pmatrix}\right],\\
\widetilde F_1(x)&=x\begin{pmatrix}0&0\\k+1&-(\alpha+\beta+3)\end{pmatrix}+\begin{pmatrix}-\frac{\beta-k+1}{\alpha+\beta-k+2}&-\frac{(\beta-k+1)(\alpha+\beta-k+1)}{\alpha+\beta-k+2}\\\frac{\alpha+1}{\alpha+\beta-k+2}&\frac{(\alpha+1)(\alpha+\beta-k+1)}{\alpha+\beta-k+2}\end{pmatrix},\\
\widetilde F_0&=\begin{pmatrix}(k+1)(\alpha+\beta-k+1)&0\\-(k+1)&0\end{pmatrix}.
\end{align*}
	Typically, in the scalar case, and for some special values of the parameters involved, the order of the differential equation satisfied by the new polynomials after a Darboux transformation is higher than 2. The remarkable fact in the matrix case is that we have a family of matrix-valued orthogonal polynomials $\widetilde P_n$ (depending on one free parameter $s_{11}$) satisfying \emph{the same second-order} differential equation with coefficients \emph{independent} of $s_{11}$. This phenomenon is not new and appeared for the first time in \cite{DdI2} using a method different than the Darboux transformation. For other examples of the bispectral property following a Darboux transformation see \cite{G3}.
\end{remark}

\end{enumerate}

\begin{remark}\label{Remp}
Once we have the explicit expression of the matrix-valued spectral measure $W$ associated with $P$ (or $\widetilde W$ associated with $\widetilde P$, or $\widehat W$ associated with $\widehat P$) we can use the Karlin-McGregor representation formula \eqref{KMcGRF} to get the $n$-step transition probability matrix $P^n$ by computing the first matrix-valued orthogonal polynomials $Q_n$. We can also compute the invariant measure $\bm\pi$ for the process $P$ using formula \eqref{InvMea}. An invariant measure for the process $P$ at hand was computed in \cite{GdI2}. Finally we can also study recurrence associated with the process $P$. According to Theorem 8.1 of \cite{GdI2}, the QBD process that results from $P$ is never positive recurrent. If $-1<\beta\leq0$, then the process is null recurrent. If $\beta>0$, then the process is transient. Therefore recurrence is independent of the value of $\alpha$. The QBD processes $\widetilde P$ and $\widehat P$ will inherit the same recurrence behavior as the original $P$, since the matrix-valued spectral measures $\widetilde W$ and $\widehat W$ will have the same behavior as $W$ at $x=1$ (see \cite{GdI2} for details).
\end{remark}

We hope that the discussion above gives an indication of the many possibilities that open up in the matrix-valued case. A relatively simple instance is discussed in the next section.

\section{An urn model for the $2\times2$ matrix-valued orthogonal polynomials}\label{SEC5}

We now give an urn model associated with the $2\times2$ matrix-valued orthogonal polynomials of Jacobi type given in the previous section. For this purpose we will focus on the simplest case of the UL block factorization $P=P_UP_L$ with block entries $X_n,Y_n,S_n,R_n$ given by \eqref{Xn2} and \eqref{Sn2}. In \cite{GPT3} one finds another urn model associated with this example, but different from the one to be given here.

From now on, it will be assumed that the parameters $\alpha,\beta$ and $k$ are nonnegative integers with $1\leq k\leq\beta$. Consider the discrete-time QBD process on $\mathbb{Z}_{\geq0}\times\{1,2\}$ whose one-step transition probability matrix $P$ is given by the coefficients $A_n,B_n$ and $C_n$ in \eqref{aan}, \eqref{bbn} and \eqref{ccn}, respectively. Consider the UL factorization $P=P_UP_L$ \eqref{Pxy} with coefficients $X_n,Y_n,S_n,R_n$ given by \eqref{Xn2} and \eqref{Sn2}. Each one of these matrices $P_U$ and $P_L$ will give rise to an experiment in terms of an urn model, which we call Experiment 1 and Experiment 2, respectively. For simplicity, we will consider these two experiments as discrete-time Markov chains on $\mathbb{Z}_{\geq0}$ with transitions between not only adjacent states but second adjacent ones too. At times $t=0,1,2,\ldots$ an urn A contains $n$ blue balls and this determines the state of our Markov chain on ${\mathbb Z}_{\geq 0}$ at that time. All the urns we use in both experiments sit in a bath consisting of an infinite number of blue and red balls.

Experiment 1 (for $P_U$) will give a discrete-time pure-birth QBD process on ${\mathbb Z}_{\geq 0}\times\{1,2\}$ with diagram given by

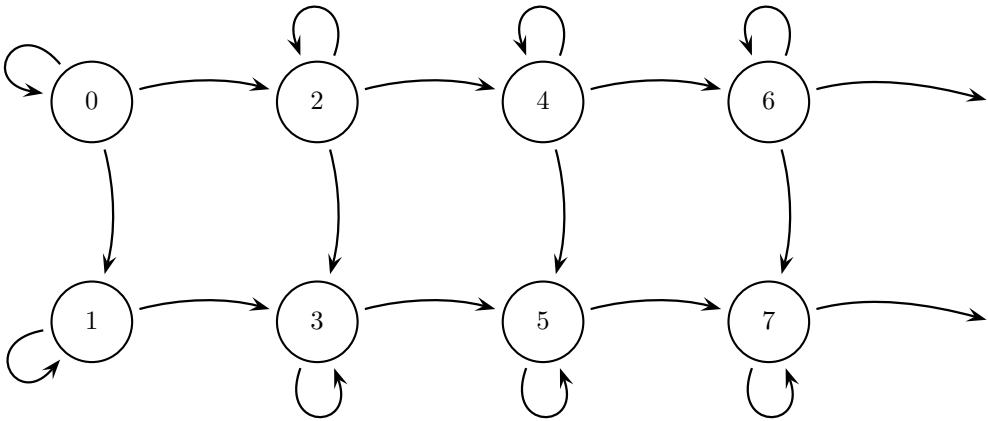
\vspace{1.0cm}
\begin{center}
$$
\begin{psmatrix}[rowsep=2.5cm,colsep=3cm]
  \cnode{.55}{0}& \cnode{.55}{2} & \cnode{.55}{4} & \cnode{.55}{6} & \pnode{8} \\
  \cnode{.55}{1} & \cnode{.55}{3} & \cnode{.55}{5} & \cnode{.55}{7} & \pnode{9} \\
\psset{nodesep=3pt,arcangle=15,labelsep=2ex,linewidth=0.3mm,arrows=->,arrowsize=1mm
3} \nccurve[angleA=130,angleB=170,ncurv=4]{0}{0}
\nccurve[angleA=190,angleB=230,ncurv=4]{1}{1}
\nccurve[angleA=70,angleB=110,ncurv=4]{2}{2}
\nccurve[angleA=70,angleB=110,ncurv=4]{4}{4}
\nccurve[angleA=70,angleB=110,ncurv=4]{6}{6}
\nccurve[angleA=250,angleB=290,ncurv=4]{3}{3}
\nccurve[angleA=250,angleB=290,ncurv=4]{5}{5}
\nccurve[angleA=250,angleB=290,ncurv=4]{7}{7} \ncarc{0}{2}
\ncarc{2}{4} \ncarc{4}{6}
\ncarc{6}{8} \ncarc{0}{1}
\ncarc{1}{3} \ncarc{2}{3}
\ncarc{3}{5}
\ncarc{4}{5}\ncarc{5}{7}
\ncarc{6}{7}\ncarc{7}{9}
\psset{labelsep=-4.25ex}\nput{90}{0}{0}
\psset{labelsep=-4.25ex}\nput{90}{2}{2}
\psset{labelsep=-4.25ex}\nput{90}{4}{4}
\psset{labelsep=-4.25ex}\nput{90}{6}{6}
\psset{labelsep=-4.25ex}\nput{90}{1}{1}
\psset{labelsep=-4.25ex}\nput{90}{3}{3}
\psset{labelsep=-4.25ex}\nput{90}{5}{5}
\psset{labelsep=-4.25ex}\nput{90}{7}{7}
\end{psmatrix}
$$
\end{center}

\vspace{-1.5cm}

This latter diagram can also be viewed as a pure-birth discrete-time Markov chain on $\mathbb{Z}_{\geq0}$ with transitions between not only adjacent states but second adjacent ones too. Let us call this chain $\{Z_t^{(1)} : t=0,1,\ldots\}$. A diagram of the same situation is now given by

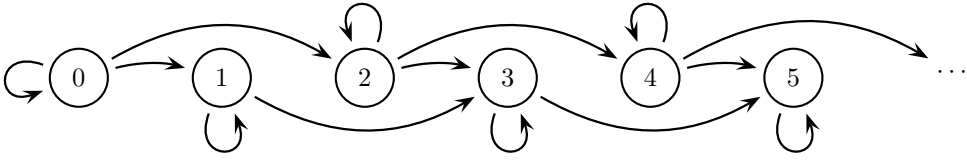
\vspace{0.5cm}

\begin{center}
$$
\begin{psmatrix}[colsep=1.9cm]
  \cnode{.4}{0}& \cnode{.4}{1} & \cnode{.4}{2}& \cnode{.4}{3}& \cnode{.4}{4}& \cnode{.4}{5} &  \rnode{6}{\Huge{\cdots}} \\
\psset{nodesep=3pt,arcangle=15,labelsep=2ex,linewidth=0.3mm,arrows=->,arrowsize=1mm
3} \nccurve[angleA=160,angleB=200,ncurv=4]{0}{0}
\nccurve[angleA=250,angleB=290,ncurv=4]{1}{1}
\nccurve[angleA=70,angleB=110,ncurv=4]{2}{2}
\nccurve[angleA=250,angleB=290,ncurv=4]{3}{3}
\nccurve[angleA=70,angleB=110,ncurv=4]{4}{4}
\nccurve[angleA=250,angleB=290,ncurv=4]{5}{5}
\nccurve[angleA=30,angleB=150,ncurv=0.8]{0}{2}
\nccurve[angleA=30,angleB=150,ncurv=0.8]{2}{4}
\nccurve[angleA=30,angleB=150,ncurv=0.8]{4}{6}
\nccurve[angleA=330,angleB=210,ncurv=0.8]{1}{3}
\nccurve[angleA=330,angleB=210,ncurv=0.8]{3}{5}
 \ncarc{0}{1}\ncarc{2}{3}\ncarc{4}{5}
\psset{labelsep=-3.40ex}\nput{90}{0}{0}
\psset{labelsep=-3.40ex}\nput{90}{1}{1}
\psset{labelsep=-3.40ex}\nput{90}{2}{2}
\psset{labelsep=-3.40ex}\nput{90}{3}{3}
\psset{labelsep=-3.40ex}\nput{90}{4}{4}
\psset{labelsep=-3.40ex}\nput{90}{5}{5}
\end{psmatrix}
$$
\end{center}

\vspace{-0.8cm}

We will construct an urn model for this last diagram. Assume the urn A contains $n$ blue balls ($n\geq0$) at time 0 (i.e. $Z_0^{(1)}=n$). The transition mechanism will depend on the parity of $n$.

\smallskip

Consider first the case where $n$ is odd and write $n=2m+1,m\geq0$. Remove $m+1$ blue balls from the urn A until we have $m$ blue balls. Take $\beta+2$ blue balls and $m+\alpha+1$ red balls from the bath and add them to the urn. Draw one ball from the urn at random with the uniform distribution. We have two possibilities:
\begin{itemize}
\item If we get a blue ball then we remove/add balls until we have $2m+3$ blue balls in urn A and start over. Therefore
$$
\mathbb{P}\left(Z_1^{(1)}=n+2 \, | \,  Z_0^{(1)}=n,\;n=2m+1\right)=\frac{m+\beta+2}{2m+\alpha+\beta+3}.
$$
Observe that this probability is given by entry $(2,2)$ of $X_m$ in \eqref{Xn2}.
\item If we get a red ball then we remove/add balls until we have $2m+1$ blue balls in urn A and start over. Therefore
$$
\mathbb{P}\left(Z_1^{(1)}=n \, | \,  Z_0^{(1)}=n,\;n=2m+1\right)=\frac{m+\alpha+1}{2m+\alpha+\beta+3}.
$$
Observe that this probability is given by entry $(2,2)$ of $Y_m$ in \eqref{Xn2}.
\end{itemize}

Consider now the case where $n$ is even and write $n=2m, m\geq0$. Again, remove $m$ blue balls from the urn A until we have $m$ blue balls. Additionally we will have two other urns, one painted in blue, which we call B, and the other one painted in red, which we call R. These urns are initially empty and will be emptied after their use in going from one time step to the next.

In urn A we add $\alpha$ blue balls and $\beta-k+1$ red balls. In urn B we place $m+\alpha+\beta-k+2$ blue balls and $m+k$ red balls. In urn R we place $m+k$ blue balls and 1 red ball. These balls are taken from the bath. Draw one ball from urn A at random with the uniform distribution. If we get a blue ball then we go to the urn B and draw a ball, while if we get a red ball then we go the urn R and draw a ball. We have three possibilities:
\begin{itemize}
\item  If we get two blue balls in a row, i.e. one from urn A and then one from urn B, then we remove/add balls until we have $2m$ blue balls in urn A and start over. Therefore
$$
\mathbb{P}\left(Z_1^{(1)}=n \, | \,  Z_0^{(1)}=n,\;n=2m\right)=\frac{(m+\alpha)(m+\alpha+\beta-k+2)}{(2m+\alpha+\beta+2)(m+\alpha+1-k+\beta)}.
$$
Observe that this probability is given by entry $(1,1)$ of $Y_m$ in \eqref{Xn2}.
\item  If we get two red balls in a row, i.e. one from urn A and then one from urn R, then we remove/add balls until we have $2m+1$ blue balls in urn A and start over. Therefore
$$
\mathbb{P}\left(Z_1^{(1)}=n+1 \, | \,  Z_0^{(1)}=n,\;n=2m\right)=\frac{\beta-k+1}{(m+\alpha+\beta-k+1)(m+k+1)}.
$$
Observe that this probability is given by entry $(1,2)$ of $Y_m$ in \eqref{Xn2}. 
\item If we get either a blue and a red ball or a red and a blue ball then we remove/add balls until we have $2m+2$ blue balls in urn A and start over. Therefore
$$
\mathbb{P}\left(Z_1^{(1)}=n+2 \, | \,  Z_0^{(1)}=n,\;n=2m\right)=\frac{(m+k)(m+\beta+2)}{(2m+\alpha+\beta+2)(m+k+1)}.
$$
Observe that this probability is given by entry $(1,1)$ of $X_m$ in \eqref{Xn2}.
\end{itemize}
We are done describing Experiment 1 and we move on to describe an unrelated experiment.

\vspace{1.0cm}

Experiment 2 (for $P_L$) will give a discrete-time pure-death QBD process on ${\mathbb Z}_{\geq 0}\times\{1,2\}$ with diagram given by

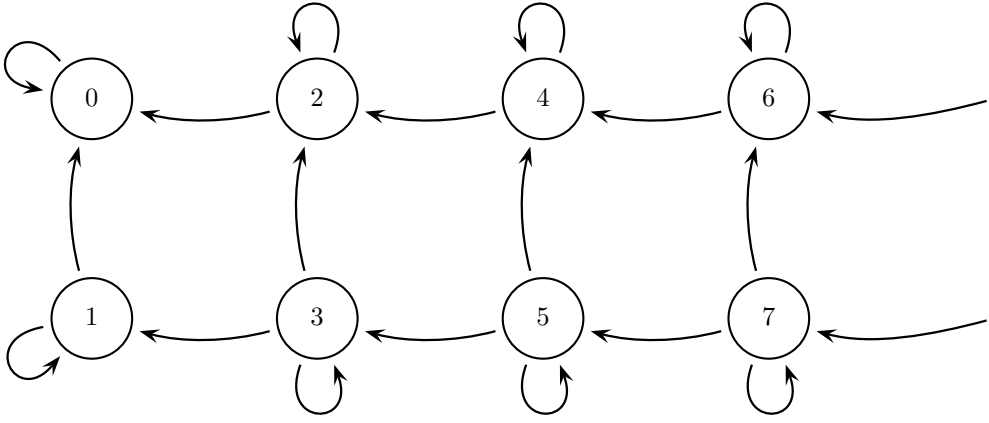
\vspace{1.0cm}

$$
\begin{psmatrix}[rowsep=2.5cm,colsep=3cm]
  \cnode{.55}{0}& \cnode{.55}{2} & \cnode{.55}{4} & \cnode{.55}{6} & \pnode{8} \\
  \cnode{.55}{1} & \cnode{.55}{3} & \cnode{.55}{5} & \cnode{.55}{7} & \pnode{9} \\
\psset{nodesep=3pt,arcangle=15,labelsep=2ex,linewidth=0.3mm,arrows=->,arrowsize=1mm
3} \nccurve[angleA=130,angleB=170,ncurv=4]{0}{0}
\nccurve[angleA=190,angleB=230,ncurv=4]{1}{1}
\nccurve[angleA=70,angleB=110,ncurv=4]{2}{2}
\nccurve[angleA=70,angleB=110,ncurv=4]{4}{4}
\nccurve[angleA=70,angleB=110,ncurv=4]{6}{6}
\nccurve[angleA=250,angleB=290,ncurv=4]{3}{3}
\nccurve[angleA=250,angleB=290,ncurv=4]{5}{5}
\nccurve[angleA=250,angleB=290,ncurv=4]{7}{7} 
\ncarc{2}{0}
\ncarc{4}{2} \ncarc{6}{4}
\ncarc{8}{6} \ncarc{1}{0}
\ncarc{3}{1} \ncarc{3}{2}
\ncarc{5}{3}
\ncarc{5}{4}\ncarc{7}{5}
\ncarc{7}{6}\ncarc{9}{7}
\psset{labelsep=-4.25ex}\nput{90}{0}{0}
\psset{labelsep=-4.25ex}\nput{90}{2}{2}
\psset{labelsep=-4.25ex}\nput{90}{4}{4}
\psset{labelsep=-4.25ex}\nput{90}{6}{6}
\psset{labelsep=-4.25ex}\nput{90}{1}{1}
\psset{labelsep=-4.25ex}\nput{90}{3}{3}
\psset{labelsep=-4.25ex}\nput{90}{5}{5}
\psset{labelsep=-4.25ex}\nput{90}{7}{7}
\end{psmatrix}
$$

\vspace{-1.4cm}

Again, this last diagram can also be viewed as a pure-death discrete-time Markov chain on $\mathbb{Z}_{\geq0}$ with transitions between not only adjacent states but second adjacent ones too, and with an absorbing state at $0$. Let us call this chain $\{Z_t^{(2)} : t=0,1,\ldots\}$. A diagram of the same situation is now given by 
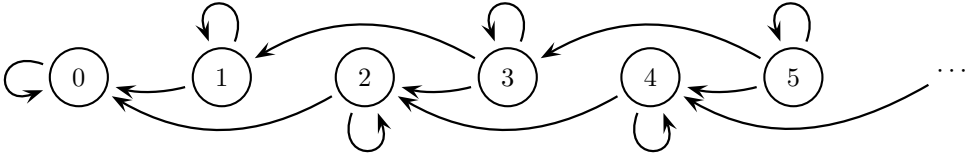
\vspace{0.7cm}
\begin{center}
$$
\begin{psmatrix}[colsep=1.9cm]
  \cnode{.4}{0}& \cnode{.4}{1} & \cnode{.4}{2}& \cnode{.4}{3}& \cnode{.4}{4}& \cnode{.4}{5} &  \rnode{6}{\Huge{\cdots}} \\
\psset{nodesep=3pt,arcangle=15,labelsep=2ex,linewidth=0.3mm,arrows=->,arrowsize=1mm
3} \nccurve[angleA=160,angleB=200,ncurv=4]{0}{0}
\nccurve[angleA=250,angleB=290,ncurv=4]{2}{2}
\nccurve[angleA=70,angleB=110,ncurv=4]{1}{1}
\nccurve[angleA=250,angleB=290,ncurv=4]{4}{4}
\nccurve[angleA=70,angleB=110,ncurv=4]{3}{3}
\nccurve[angleA=70,angleB=110,ncurv=4]{5}{5}
\nccurve[angleA=150,angleB=30,ncurv=0.8]{3}{1}
\nccurve[angleA=150,angleB=30,ncurv=0.8]{5}{3}
\nccurve[angleA=210,angleB=330,ncurv=0.8]{6}{4}
\nccurve[angleA=210,angleB=330,ncurv=0.8]{2}{0}
\nccurve[angleA=210,angleB=330,ncurv=0.8]{4}{2}
 \ncarc{1}{0}\ncarc{3}{2}\ncarc{5}{4}
\psset{labelsep=-3.40ex}\nput{90}{0}{0}
\psset{labelsep=-3.40ex}\nput{90}{1}{1}
\psset{labelsep=-3.40ex}\nput{90}{2}{2}
\psset{labelsep=-3.40ex}\nput{90}{3}{3}
\psset{labelsep=-3.40ex}\nput{90}{4}{4}
\psset{labelsep=-3.40ex}\nput{90}{5}{5}
\end{psmatrix}
$$
\end{center}

\vspace{-0.6cm}

We will construct an urn model for this last diagram. Assume that urn A contains $n$ blue balls ($n\geq1$) at time 0 (i.e. $Z_0^{(0)}=n$). The state $n=0$ is absorbing. Consider first the case where $n$ is even and write $n=2m, m\geq1$. Remove $m$ blue balls from the urn A until we have $m$ blue balls. Take $m+\alpha+\beta+1$ red balls from the bath and add them to the urn. Draw one ball from the urn at random with the uniform distribution. We have two possibilities:
\begin{itemize}
\item If we get a blue ball then we remove/add balls until we have $2m-2$ blue balls in urn A and start over. Therefore
$$
\mathbb{P}\left(Z_1^{(2)}=n-2 \, | \,  Z_0^{(2)}=n,\;n=2m\right)=\frac{m}{2m+\alpha+\beta+1}.
$$
Observe that this probability is given by entry $(1,1)$ of $R_m$ in \eqref{Sn2}.
\item If we get a red ball then we remove/add balls until we have $2m$ blue balls at the urn A and start over. Therefore
$$
\mathbb{P}\left(Z_1^{(2)}=n \, | \,  Z_0^{(2)}=n,\;n=2m\right)=\frac{m+\alpha+\beta+1}{2m+\alpha+\beta+1}.
$$
Observe that this probability is given by entry $(1,1)$ of $S_m$ in \eqref{Sn2}.
\end{itemize}

Consider now the case where $n$ is odd and write $n=2m+1, m\geq0$. Again, remove $m+1$ blue balls from the urn A until we have $m$ blue balls. Additionally we will have two other urns, one painted in blue, which we call B, and the other one painted in red, which we call R. Again, these urns are initially empty and will be emptied after their use in going from one time step to the next.

In urn A we add $m+\alpha+\beta-k+1$ blue balls and $1$ red ball. In urn B we place $m+\alpha+\beta+2$ blue balls and $m$ red balls. In urn R we place $m$ blue balls and $k$ red balls. Draw one ball from urn A at random with the uniform distribution. If we get a blue ball then we go to the urn B and draw a ball, while if we get a red ball then we go the urn R and draw a ball. We have three possibilities:
\begin{itemize}
\item If we get two blue balls in a row, i.e. one from urn A and then one from urn B, then we remove/add balls until we have $2m+1$ blue balls in urn A and start over. Therefore
$$
\mathbb{P}\left(Z_1^{(2)}=n \, | \,  Z_0^{(2)}=n,\;n=2m+1\right)=\frac{(m+\alpha+\beta+2)(m+\alpha+1-k+\beta)}{(2m+\alpha+\beta+2)(m+\alpha+\beta-k+2)}.
$$
Observe that this probability is given by entry $(2,2)$ of $S_m$ in \eqref{Sn2}.
\item  If we get two red balls in a row, i.e. one from urn A and then one from urn R, then we remove/add balls until we have $2m$ blue balls in urn A and start over. Therefore
$$
\mathbb{P}\left(Z_1^{(2)}=n-1 \, | \,  Z_0^{(2)}=n,\;n=2m+1\right)=\frac{k}{(m+\alpha+\beta-k+2)(m+k)}.
$$
Observe that this probability is given by entry $(2,1)$ of $S_m$ in \eqref{Sn2}.
\item If we get either a blue and a red ball or a red and a blue ball then we remove/add balls until we have $2m-1$ blue balls in urn A and start over. Therefore
$$
\mathbb{P}\left(Z_1^{(2)}=n-2 \, | \,  Z_0^{(2)}=n,\;n=2m+1\right)=\frac{m(m+k+1)}{(2m+\alpha+\beta+2)(m+k)}.
$$
Observe that this probability is given by entry $(2,2)$ of $R_m$ in \eqref{Sn2}.
\end{itemize}

The urn model for $P$ (on $\mathbb{Z}_{\geq0}$) will be the composition of Experiment 1 and then Experiment 2, while the urn model for the Darboux transformation $\widetilde P$ \eqref{DarbT} proceeds in the reversed order. Observe from Remark \ref{Remp} and since $\alpha,\beta$ and $k$ are nonnegative integers with $1\leq k\leq\beta$ that this urn model will always be transient. Similar urn models can be derived for the LU factorization with small modifications.


\begin{thebibliography}{99}



\bibitem{AW} Askey, R. and Wilson, J., {\em Some basic hypergeometric orthogonal polynomials that generalize Jacobi polynomials},  Mem. Amer. Math. Soc. {\bf 54} (1985), no. 319.




\bibitem{DRSZ} Dette, H., Reuther, B., Studden, W. and Zygmunt, M., {\em Matrix measures and random walks with a block tridiagonal transition matrix}, SIAM J. Matrix Anal. Applic. \textbf{29}, No. 1 (2006), 117--142.

\bibitem{DdI2} Durán, A. J. and de la Iglesia, M. D., 
\textit{Second order differential operators having several families of orthogonal matrix polynomials as eigenfunctions}, Internat. Math. Research Notices, Vol. \textbf{2008}, Article ID rnn084, 24 pages.






\bibitem{E} Ehrenfest, P. and Eherenfest, T., {\em \"Uber zwei bekannte Einw\"ande gegen das
Boltzmannsche H-Theorem}, Physikalische Zeitschrift, vol.~8 (1907), 311--314.


\bibitem{F}  Feller, W., {\em An introduction to Probability Theory and its Applications}, vol.~1, 3rd edition, Wiley, 1967.



\bibitem{G2}\textrm{Grünbaum, F. A.},
\textit{Random walks and orthogonal polynomials: some challenges}, Probability, Geometry and Integrable Systems, MSRI Publication, volumen \textbf{55}, 2007.


\bibitem{G4} \textrm{Gr\"unbaum, F. A.}, \emph{An urn model associated with Jacobi polynomials}, Commun. Applied Math. Comput. Sciences \textbf{5} (2010), no. 1, 55--63.

\bibitem{G3}\textrm{Grünbaum, F. A.}, \emph{The Darboux process and a noncommutative bispectral problem: some explorations and challenges}, in E.P. van den Ban and J.A.C. Kolk (eds.), \emph{Geometric Aspects of Analysis and Mechanics: In Honor of the 65th Birthday of Hans Duistermaat}, Progress in Mathematics 292, Springer, 2011.


\bibitem{GdI2} Gr\"unbaum, F. A. and de la Iglesia, M. D., \textit{Matrix valued orthogonal polynomials arising from group representation theory and a family of quasi-birth-and-death processes}, SIAM J. Matrix Anal. Applic. \textbf{30}, No. 2 (2008), 741--761.

\bibitem{GdI3} Gr\"unbaum, F. A. and de la Iglesia, M. D.,
\textit{Stochastic LU factorizations, Darboux transformations and urn models}, submitted, see arXiv:1706.02617.

\bibitem{GPT1} Grünbaum, F. A., Pacharoni, I. and Tirao, J. A.,
{\em Matrix valued spherical functions associated to the complex
projective plane}, J. Functional Analysis {\bf 188} (2002),
350--441.

\bibitem{GPT2} Grünbaum, F. A., Pacharoni, I. and Tirao, J. A.,
{\em A matrix valued solution to Bochner's problem}, J. Physics A:
Math. Gen. {\bf 34} (2001), 10647--10656.


\bibitem{GPT3} Gr\"unbaum, F. A., Pacharoni, I. and Tirao, J.
A., {\em Two stochastic models of a random walk in the U($n$)-spherical duals of U($n+1$)}, Ann. Mat. Pura Appl. \textbf{192} (2013), no. 3, 447--473. 

\bibitem{dI2} de la Iglesia, M. D., \emph{A note on the invariant distribution of a quasi-birth-and-death process}, J. Phys. A: Math. Theor. \textbf{44} (2011), 135201 (9pp).


\bibitem{ILMV}  Ismail, M.E.H., Letessier, J., Masson, D. and Valent, G., {\em Birth and death processes and orthogonal polynomials}, in Orthogonal Polynomials, P.~Nevai (editor),  Kluwer Acad. Publishers (1990), 229--255.

\bibitem{KMc2}  Karlin, S. and McGregor, J., {\em The differential equations of birth and death processes, and the Stieltjes moment problem}, Trans. Amer. Math. Soc., \textbf{85} (1957), 489--546.

\bibitem{KMc3}  Karlin, S. and McGregor, J., {\em The classification of birth-and-death processes}, Trans. Amer. Math. Soc., \textbf{86} (1957), 366--400.


\bibitem{KMc6}  Karlin, S. and McGregor, J., {\em Random walks}, IIlinois J. Math., \textbf{3} (1959), 66--81.


\bibitem{K2} Kre$\breve{\mbox{{\i}}}$n, M. G., {\em Infinite $J$-matrices and a matrix moment problem}, Dokl. Akad. Nauk SSSR {\bf 69} No. 2 (1949), 125--128.

\bibitem{K1}  Kre$\breve{\mbox{{\i}}}$n, M. G., {\em  Fundamental aspects of the representation theory of hermitian operators with deficiency index $(m,m)$},  AMS Translations, Series 2, \textbf{97} (1971), Providence, Rhode Island, 75--143.

\bibitem{LaR} Latouche, G. and Ramaswami, V., {\em Introduction to Matrix Analytic Methods in Stochastic Modeling}, ASA-SIAM Series on Statistics and Applied Probability, 1999.

\bibitem{MS} Matveev, V.B. and Salle, M.A., \emph{Differential-difference evolution equations II: Darboux transformation for the Toda lattice}, Lett. Math. Phys. \textbf{3} (1979) 425--429.


\bibitem{Neu} Neuts, M. F., {\em Structured Stochastic Matrices of $M/G/1$ Type and Their Applications}, Marcel Dekker, New York, 1989.

\bibitem{PT1}  Pacharoni, I. and  Tirao, J. A.,
{\em Matrix-valued orthogonal polynomials arising from the complex
projective space}. Constr. Approx.  \textbf{25} (2007), pp.
177--192.


\end{thebibliography}
\end{document}